\documentclass[a4paper,11pt,twoside]{amsart}

\usepackage[T1]{fontenc}
\usepackage[latin1]{inputenc}% Codage du fichier en ISO-Latin-1 (accents...)
\usepackage[francais]{babel}% Utilisation du Franà¯Â¿Â½is comme langue principale, de l'anglais comme langue secondaire
\usepackage{xspace} % Espaces insécables avant les ponctuations doubles
\usepackage{amsfonts}
\usepackage{amsmath,amssymb,amsthm,,amscd}% Formules mathématiques
\usepackage[all]{xy}
\usepackage{graphicx}

\title[Quotients compacts des groupes ultramétriques de rang un]{Quotients compacts des\\ groupes ultramétriques de rang un}
\author{Fanny Kassel}

\theoremstyle{plain}
\newtheorem{prop}[table]{Proposition}
\newtheorem{theo}[table]{Th\'eor\`eme}
\newtheorem{coro}[table]{Corollaire}
\newtheorem{lem}[table]{Lemme}
\newtheorem{rema}[table]{Remarque}
\newenvironment{dem}{\noindent{\textbf{D\'emonstration.}}}{\hfill \qedsymbol}

\newcommand{\C}{\mathbb{C}}
\newcommand{\R}{\mathbb{R}}
\newcommand{\Q}{\mathbb{Q}}
\newcommand{\Z}{\mathbb{Z}}
\newcommand{\N}{\mathbb{N}}
\newcommand{\kkk}{\mathbf{k}}
\newcommand{\F}{\mathbb{F}}
\newcommand{\SL}{\mathrm{SL}}
\newcommand{\PSL}{\mathrm{PSL}}

\newcommand{\Ad}{\mathrm{Ad}}
\newcommand{\g}{\mathfrak{g}}

\newcommand{\dd}{\underline{d}}
\newcommand{\Hom}{\mathrm{Hom}}

\newcommand{\Isom}{\mathrm{Isom}}
\newcommand{\D}{\mathcal{D}}
\newcommand{\pr}{\mathrm{pr}}
\newcommand{\OS}{\mathrm{OS}_n}

\begin{document}
\maketitle
\numberwithin{equation}{section}
\numberwithin{table}{section}

%%%%%%%%%%%%%%%%%%%%%%%%%%%%%%%%%%%%%%%%%%%%%%%%%%%
%%%%%%%%%%%%%%%%%%%%%%%%%%%%%%%%%%%%%%%%%%%%%%%%%%%
\begin{abstract}
Soit $G$ l'ensemble des $\kkk$-points d'un groupe algébrique semi-simple connexe de $\kkk$-rang~un sur un corps local ultramétrique~$\kkk$.
Nous décrivons tous les sous-groupes discrets de type fini sans torsion de~$G\times\nolinebreak G$ qui agissent proprement et cocompactement sur~$G$ par multiplication à gauche et à droite.
Nous montrons qu'après une petite déformation dans~$G\times G$ un tel sous-groupe discret agit encore librement, proprement et cocompactement sur~$G$.

\bigskip

\noindent
\textsc{Abstract.}
Let $G$ be the set of $\kkk$-points of a connected semisimple algebraic group of $\kkk$-rank~one over a nonarchimedean local field~$\kkk$.
We describe all finitely generated torsion-free discrete subgroups of~$G\times\nolinebreak G$ acting properly discontinuously and cocompactly on~$G$ by left and right multiplication.
We prove that after a small deformation in~$G\times G$ such a discrete subgroup keeps acting freely, properly discontinuously, and cocompactly on~$G$.
\end{abstract}

%%%%%%%%%%%%%%%%%%%%%%%%%%%%%%%%%%%%%%%%%%%%%%%%%%%
%%%%%%%%%%%%%%%%%%%%%%%%%%%%%%%%%%%%%%%%%%%%%%%%%%%
\section{Introduction}\label{Introduction}

Soient~$\kkk$ un corps local et $G$ l'ensemble des $\kkk$-points d'un groupe algébrique semi-simple connexe~$\mathbf{G}$ sur~$\kkk$, de $\kkk$-rang~un.
Nous nous intéressons aux sous-groupes discrets~$\Gamma$ de~$G\times G$ qui agissent librement, proprement et cocompactement sur~$G$ par multiplication à gauche et à droite.
De manière équivalente, ce sont les sous-groupes discrets~$\Gamma$ de~$G\times G$ qui agissent librement, proprement et cocompactement sur l'espace homogène~$(G\times G)/\Delta_G$, où $\Delta_G$ désigne la diagonale de~$G\times G$.
Pour un tel groupe~$\Gamma$, on dit que le quotient $\Gamma\backslash (G\times G)/\Delta_G$ est une \textit{forme de Clifford-Klein compacte}, ou plus simplement un \textit{quotient compact}, de~$(G\times G)/\Delta_G$.

Pour $G=\PSL_2(\R)$, par exemple, les quotients compacts de~$(G\times G)/\Delta_G$ ont été largement étudiés par W.~M.~Goldman \cite{gol}, T.~Kobayashi \cite{ko2}, \cite{ko3}, R.~S.~Kulkarni et F.~Raymond \cite{kr}, F.~Salein \cite{sa1}, \cite{sa2}.
Ils apparaissent naturellement en géométrie : les variétés anti-de Sitter compactes de dimension~3, c'est-à-dire les variétés lorentziennes compactes de dimension~3 de courbure sectionnelle constante égale à~$-1$, sont exactement les revêtements finis des quotients compacts de $(\PSL_2(\R)\times\PSL_2(\R))/\Delta_{\PSL_2(\R)}$.
Cela résulte de la complétude de ces variétés~\cite{kli} et d'un résultat de finitude du niveau~\cite{kr} (\textit{cf.} l'introduction de~\cite{sa2}).
Pour $G=\SL_2(\C)$, les quotients compacts de~$(G\times G)/\Delta_G$ ont également été étudiés par \'E.~Ghys \cite{ghy} en lien avec les déformations de structures complexes sur les variétés compactes homogènes sous~$G$.

Dans cet article, nous nous intéressons au cas d'un corps local~$\kkk$ ultramétrique, c'est-à-dire d'une extension finie de~$\Q_p$ ou du corps~$\F_q((t))$ des séries de Laurent formelles à coefficients dans un corps fini~$\F_q$.
L'une de nos motivations vient de l'étude des quotients compacts de la quadrique
$$\big\{ (x_1,x_2,x_3,x_4)\in\kkk^4,\quad x_1^2 - x_2^2 + x_3^2 - x_4^2 = 1\big\} ,$$
qui s'identifie à l'espace homogène $(\SL_2(\kkk)\times\SL_2(\kkk))/\Delta_{\SL_2(\kkk)}$ (\textit{cf.} \cite{kas},~\S~5.3).

Pour tout groupe algébrique semi-simple connexe~$\mathbf{G}$ sur~$\kkk$, de $\kkk$-rang~un, l'espace homogène $(G\times G)/\Delta_G$ admet un quotient compact.
En effet, $G$ admet un réseau cocompact sans torsion~$\Gamma_0$ (\cite{lub}, th.~A) ; le groupe $\Gamma_0\times\nolinebreak\{ 1\} $ agit librement, proprement et cocompactement sur $(G\times G)/\Delta_G$.
Les questions suivantes se posent alors naturellement :
\begin{enumerate}
	\item décrire tous les quotients compacts~$\Gamma\backslash (G\times G)/\Delta_G$ ;
	\item comprendre le comportement de ces quotients compacts lorsque l'on déforme~$\Gamma$ dans~$G$; en particulier, déterminer si l'action de~$\Gamma$ reste propre et cocompacte ;
	\item établir l'existence d'un quotient compact $\Gamma\backslash (G\times G)/\Delta_G$ tel que $\Gamma$~soit Zariski-dense dans~$G\times G$.
\end{enumerate}
Dans cet article, nous répondons aux trois questions.
Nous obtenons ainsi des analogues ultramétriques de résultats connus pour~$\PSL_2(\R)$.

%%%%%%%%%%%%%%%%%%%%%%%%%
\subsection{\'Enoncé des résultats principaux}

Rappelons que le groupe~$G$ admet une décomposition de Cartan de la forme~$G=KZ^+K$, où~$K$ est un sous-groupe compact maximal de~$G$ et~$Z^+$ une chambre de Weyl de l'ensemble des $\kkk$-points du centralisateur d'un $\kkk$-tore $\kkk$-déployé maximal de~$\mathbf{G}$ ; à cette décomposition de Cartan est naturellement associée une projection de Cartan $\mu : G\rightarrow\R^+$ (\textit{cf.} paragraphe \ref{Cartan}).
Par exemple, si~$G=\SL_2(\kkk)$, le théorème de la base adaptée induit la décomposition de Cartan~$G=KZ^+K$, où~$K=\SL_2(\mathcal{O})$ est l'ensemble des matrices de déterminant~un à coefficients dans l'anneau des entiers de~$\kkk$ et $Z^+$ l'ensemble des matrices diagonales~$\mathrm{diag}(a,a^{-1})\in\SL_2(\kkk)$ telles que~$a$ soit de valeur absolue~$\geq 1$ ; la projection de Cartan associée~$\mu : G\rightarrow\R^+$ envoie la matrice~$\mathrm{diag}(a,a^{-1})$ sur~$2\,|\omega(a)|$, où~$\omega$ désigne une valuation (additive) fixée de~$\kkk$.

On sait décrire en fonction de~$\mu$ tous les sous-groupes discrets sans torsion~$\Gamma$ de~$G\times G$ qui agissent proprement sur~$(G\times G)/\Delta_G$ (\cite{kas}, th.~1.3).
Dans cet article, nous établissons un critère pour que le quotient $\Gamma\backslash (G\times G)/\Delta_G$ soit compact (théorème~\ref{equivalence cocompacite}).
Nous obtenons ainsi une description de tous les quotients compacts de~$(G\times G)/\Delta_G$ par un groupe~$\Gamma$ discret de type fini sans torsion.

\begin{theo}\label{theoreme quotients compacts}
Soient $\kkk$ un corps local ultramétrique, $G$ l'ensemble des $\kkk$-points d'un $\kkk$-groupe algébrique semi-simple connexe de $\kkk$-rang~un, $\Delta_G$~la diagonale de~$G\times G$ et $\mu : G\rightarrow\R^+$ une projection de Cartan de~$G$.\linebreak
\`A la permutation près des deux facteurs de~$G\times G$, les sous-groupes discrets de type fini sans torsion de~$G\times G$ agissant proprement et cocompactement sur~$(G\times G)/\Delta_G$ sont exactement les graphes de la forme
$$\Gamma = \big\{ (\gamma,\rho(\gamma)),\ \gamma\in\Gamma_0\big\} ,$$
où $\Gamma_0$ est un réseau cocompact sans torsion de~$G$ et $\rho : \Gamma_0\rightarrow G$ un morphisme de groupes qui est admissible, au sens où pour tout $R>0$ on a $\mu(\rho(\gamma))\leq\nolinebreak\mu(\gamma)-R$ pour presque tout $\gamma\in\Gamma_0$.
\end{theo}

On dit ici qu'une propriété est vraie pour \textit{presque tout}~$\gamma\in\Gamma_0$ si elle est vraie pour tout~$\gamma\in\Gamma_0$ en dehors d'un ensemble fini.

Cette description est spécifique au rang~un : au paragraphe~\ref{Rang superieur} nous donnons en caractéristique nulle un exemple de quotient compact $\Gamma\backslash (G\times\nolinebreak G)/\Delta_G$ où $\mathrm{rang}_{\kkk}(\mathbf{G})\geq 2$ et où $\Gamma$ est le produit de deux sous-groupes infinis de~$G$.

Le théorème \ref{theoreme quotients compacts} est également valable pour $\kkk=\R$ d'après \cite{ko2}, th.~2, et \cite{ko1}, cor.~5.5, ainsi que \cite{kas}, th.~1.3.
Les arguments de~\cite{ko1} utilisent la dimension cohomologique de~$\Gamma$.
Dans le cas ultramétrique, ces arguments ne conviennent pas ; nous les remplaçons par des raisonnements géométriques sur l'arbre de Bruhat-Tits de~$G$.

Nous apportons ensuite une réponse positive à la question~(2).

\begin{theo}\label{ouvert}
Soient $\kkk$ un corps local ultramétrique, $G$ l'ensemble des \linebreak $\kkk$-points d'un $\kkk$-groupe algébrique semi-simple connexe de $\kkk$-rang~un et $\Gamma$ un sous-groupe discret de type fini sans torsion de~$G\times G$.
Si $\Gamma$ agit proprement et cocompactement sur~$(G\times G)/\Delta_G$, alors il existe un voisinage~$\mathcal{U}$ de l'inclusion canonique dans~$\Hom(\Gamma,G\times G)$ tel que pour tout~$\varphi\in\mathcal{U}$, le groupe~$\varphi(\Gamma)$ agisse librement, proprement et cocompactement sur $(G\times G)/\Delta_G$.
\end{theo}

On note ici~$\Hom(\Gamma,G\times G)$ l'ensemble des morphismes de groupes de~$\Gamma$ dans~$G\times G$, muni de la topologie compacte-ouverte.
Pour toute partie génératrice finie~$F$ de~$\Gamma$, une suite~$(\varphi_n)\in\Hom(\Gamma,G\times G)^{\N}$ converge vers un élément~$\varphi\in\Hom(\Gamma,G\times G)$ si et seulement si $\varphi_n(\gamma)\rightarrow\varphi(\gamma)$ pour tout~$\gamma\in F$.

Le théorème~\ref{ouvert} est également valable pour~$\kkk=\R$ et~$G=\PSL_2(\R)$ par la complétude des variétés anti-de Sitter compactes \cite{kli} et par un principe, dû à Ehresmann, de déformation des holonomies de $(G,X)$-structures sur les variétés compactes (\textit{cf.} \cite{sa1}).

Notons qu'aucun analogue du théorème~\ref{ouvert} n'est connu pour les autres groupes réels semi-simples de rang~un, même pour~$G=\SL_2(\C)$.
Le seul résultat connu pour ces groupes est l'existence d'un voisinage du morphisme trivial dans $\Hom(\Gamma_0,G)$ formé de morphismes admissibles (\cite{ko3}, th.~2.4).

Enfin nous répondons positivement à la question~(3).

\begin{prop}\label{coro Zariski-dense}
Soient $\kkk$ un corps local ultramétrique, $G$ l'ensemble des $\kkk$-points d'un $\kkk$-groupe algébrique semi-simple connexe de $\kkk$-rang~un et $\Delta_G$ la diagonale de $G\times G$.
Il existe un sous-groupe discret~$\Gamma$ de~$G\times G$ qui agit librement, proprement et cocompactement sur~$(G\times G)/\Delta_G$ et qui est Zariski-dense dans~$G\times G$.
On peut choisir~$\Gamma$ de sorte qu'aucune de ses deux projections naturelles sur~$G$ ne soit bornée.
\end{prop}

%%%%%%%%%%%%%%%%%%%%%%%%%
\subsection{Stratégie de démonstration}\label{Strategie de démonstration}

Soient~$\kkk$ un corps local ultramétrique, $G$ l'ensemble des $\kkk$-points d'un $\kkk$-groupe algébrique semi-simple connexe de $\kkk$-rang~un et $\Gamma_0$ un réseau cocompact sans torsion de~$G$.
Pour démontrer les théorèmes~\ref{theoreme quotients compacts} et~\ref{ouvert}, nous étudions l'action de~$G$ sur son arbre de Bruhat-Tits, qui est un arbre simplicial sur lequel~$G$ agit proprement par isométries, avec un domaine fondamental compact.
Nous rappelons la construction et les principales propriétés de cet arbre au paragraphe~\ref{Cartan}.

Soit~$G=KZ^+K$ une décomposition de Cartan de~$G$.
Pour démontrer le théorème~\ref{theoreme quotients compacts}, nous associons à tout élément~$g\in G\smallsetminus K$ un point~$\zeta_g^-$ du bord de l'arbre de Bruhat-Tits de~$G$, obtenu à partir d'une décomposition de Cartan de~$g$.
L'étude des points~$\zeta_g^-$ nous permet de montrer que sous les hypothèses du théorème~\ref{theoreme quotients compacts}, si $\Gamma_0$ est un sous-groupe discret de type fini sans torsion de~$G$ et $\rho : \Gamma_0\rightarrow G$ un morphisme de groupes admissible, et si le graphe
$$\Gamma = \big\{ (\gamma,\rho(\gamma)),\ \gamma\in\Gamma_0\big\} $$
agit proprement et cocompactement sur~$(G\times G)/\Delta_G$, alors $\Gamma_0$ est un réseau cocompact de~$G$.
Il s'agit du point le plus délicat de la démonstration du théorème~\ref{theoreme quotients compacts}, cette démonstration faisant l'objet de la partie~\ref{Cocompacite}.

Quant au théorème~\ref{ouvert}, il découle d'un résultat plus général sur les groupes d'isométries d'arbres réels simpliciaux (proposition~\ref{max sur un ensemble fini}).
Plus précisément, pour tout arbre réel simplicial~$X$ et toute isométrie~$g$ de~$X$, notons
\begin{equation}\label{definition lambda}
\lambda(g) = \inf_{x\in X} d(x,g\cdot x) \geq 0
\end{equation}
la \textit{longueur de translation} de~$g$, où~$d$ désigne la distance de~$X$.
L'application $\lambda : \Isom(X)\rightarrow\R^+$ ainsi définie est continue.
Par exemple, si~$X$ est l'arbre de Bruhat-Tits de~$G=\SL_2(\kkk)$, on a $\lambda(g)=|\omega(a_g)-\omega(a'_g)|$ pour tout~$g\in\nolinebreak G$, où~$a_g$ et~$a'_g$ désignent les deux valeurs propres de~$g$ et~$\omega$ une valuation (additive) fixée sur~$\kkk$.
Pour tout arbre réel simplicial~$X$ et tout sous-groupe discret sans torsion~$\Gamma_0$ de~$\Isom(X)$, on a~$\lambda(\gamma)>0$ pour tout~$\gamma\in\Gamma_0\smallsetminus\{ 1\} $.
Dans la partie~\ref{Applications equivariantes entre arbres}, nous montrons que si~$\Gamma_0$ agit cocompactement sur~$X$, alors pour tout arbre réel simplicial~$X'$ et tout morphisme de groupes~$\rho : \Gamma_0\rightarrow\Isom(X')$, la borne supérieure des quotients $\lambda(\rho(\gamma))/\lambda(\gamma)$, où~$\gamma\in\Gamma_0\smallsetminus\{ 1\} $, est égale à la plus petite constante de Lipschitz d'une application continue $f : X\rightarrow X'$ qui est $\rho$\textit{-équivariante}, au sens où
$$f(\gamma\cdot x) = \rho(\gamma)\cdot f(x)$$
pour tout~$\gamma\in\Gamma_0$ et tout~$x\in X$.
Nous montrons que cette borne supérieure est atteinte sur une partie finie~$F$ de~$\Gamma_0\smallsetminus\{ 1\} $ indépendante de~$X'$ et de~$\rho$, ce qui généralise un résultat non publié de T.~White sur l'outre-espace, dont une preuve a été donnée par S.~Francaviglia et A.~Martino~\cite{fm}.
Nous en déduisons le résultat suivant.

\begin{theo}\label{condition d'admissibilite}
Soient $\kkk$ un corps local ultramétrique, $G$ l'ensemble des $\kkk$-points d'un $\kkk$-groupe algébrique semi-simple connexe de $\kkk$-rang~un et $\Gamma_0$ un réseau cocompact sans torsion de~$G$.
Pour tout morphisme de groupes $\rho : \Gamma_0\rightarrow G$, notons~$C_{\rho}\geq 0$ la borne inférieure des réels~$t\geq 0$ pour lesquels il existe~$t'\geq 0$ tel que $\mu(\rho(\gamma))\leq t\mu(\gamma)+t'$ pour tout~$\gamma\in\Gamma_0$.
\begin{enumerate}
	\item Soit $F$ le sous-ensemble~fini de~$\Gamma_0\smallsetminus\{ 1\} $ formé des éléments~$\gamma$ tels que $\mu(\gamma)\leq 4N$, où $N$ désigne le nombre de sommets de~$\Gamma_0\backslash X$.
	Pour tout morphisme $\rho : \Gamma_0\rightarrow G$ on a
	$$C_{\rho}\ =\ \sup_{\gamma\in\Gamma_0\smallsetminus\{ 1\} } \frac{\lambda(\rho(\gamma))}{\lambda(\gamma)}\ =\ \max_{\gamma\in F}\, \frac{\lambda(\rho(\gamma))}{\lambda(\gamma)},$$
	où~$\lambda : G\rightarrow\R^+$ désigne l'application définie par~\emph{(\ref{definition lambda})} pour l'action de~$G$ sur son arbre de Bruhat-Tits.
	\item Un morphisme $\rho : \Gamma_0\rightarrow G$ est admissible si et seulement si~$C_{\rho}<1$.
\end{enumerate}
\end{theo}

L'ensemble~$F$ ci-dessus est bien fini car l'application~$\mu$ est propre et le groupe~$\Gamma_0$ discret dans~$G$.
Ainsi, le caractère admissible d'un morphisme $\rho : \Gamma_0\rightarrow G$ dépend d'un nombre \textit{fini} de conditions \textit{ouvertes}, ce qui prouve que l'ensemble des morphismes admissibles est ouvert dans~$\Hom(\Gamma_0,G)$ pour la topologie compacte-ouverte.
Ceci implique le théorème~\ref{ouvert} en utilisant une propriété des déformations de réseaux cocompacts sans torsion de~$G$ (lemme~\ref{Reste un reseau cocompact sans torsion}).

La condition suffisante du point~(2) est immédiate sachant que $\mu$ est propre et $\Gamma_0$ discret : le morphisme $\rho$ est admissible dès que~$C_{\rho}<1$.
De même, si~$\rho$ est admissible on a l'inégalité large~$C_{\rho}\leq 1$.
Toute la difficulté du point~(2) tient au fait que nous voulons obtenir une inégalité stricte.

%%%%%%%%%%%%%%%%%%%%%%%%%
\subsection{Un résultat complémentaire}

Une autre application du théorème~\ref{condition d'admissibilite} est la suivante.

\begin{coro}\label{fidele et discrete}
Soient $\kkk$ un corps local ultramétrique, $G$ l'ensemble des \linebreak $\kkk$-points d'un $\kkk$-groupe algébrique semi-simple connexe de $\kkk$-rang~un et $\Gamma_0$ un réseau cocompact sans torsion de~$G$.
Il n'existe pas de morphisme de groupes admissible~$\rho : \Gamma_0\rightarrow G$ qui soit injectif d'image discrète.
\end{coro}

Le corollaire~\ref{fidele et discrete} est également valable pour~$\kkk=\R$ et~$G=\PSL_2(\R)$ : cela résulte de l'existence, due à W.~Thurston~\cite{thu}, d'une certaine ``distance asymétrique'' sur l'espace de Teichmüller de~$\Gamma_0\backslash\mathbb{H}$, où~$\mathbb{H}$ désigne le demi-plan de Poincaré (\textit{cf.} \cite{sa2}, \S~4.1).

Pour démontrer le corollaire \ref{fidele et discrete} pour un corps local~$\kkk$ ultramétrique, nous remplaçons l'espace de Teichmüller de~$\Gamma_0\backslash\mathbb{H}$ par l'outre-espace de même rang que le groupe libre~$\Gamma_0$ (\textit{cf.} partie~\ref{outre-espace}).

%%%%%%%%%%%%%%%%%%%%%%%%%
\subsection{Plan de l'article}

La partie~\ref{Rappels} est consacrée à quelques rappels sur les décompositions de Cartan et l'arbre de Bruhat-Tits de~$G$ d'une part, sur les isométries d'arbres réels simpliciaux d'autre part.
Dans la partie~\ref{Cocompacite}, nous démontrons le théorème~\ref{theoreme quotients compacts} en étudiant les points~$\zeta_g^-$ du bord de l'arbre de Bruhat-Tits de~$G$ mentionnés précédemment.
La partie~\ref{Applications equivariantes entre arbres} établit un résultat général sur les groupes d'isométries d'arbres réels simpliciaux (proposition~\ref{max sur un ensemble fini}), dont nous déduisons les théorèmes~\ref{condition d'admissibilite} puis~\ref{ouvert} dans la partie~\ref{Deformation des quotients compacts}.
Dans la partie~\ref{outre-espace} nous décrivons le lien entre la proposition~\ref{max sur un ensemble fini} et l'outre-espace et démontrons le corollaire~\ref{fidele et discrete}.
Enfin, la partie~\ref{Zariski-dense} est consacrée à la démonstration de la proposition~\ref{coro Zariski-dense}.

%%%%%%%%%%%%%%%%%%%%%%%%%
\subsection*{Remerciements}

Je remercie vivement Yves Benoist pour de nombreuses discussions, ainsi que Frédéric Paulin et Mladen Bestvina pour leurs indications sur l'outre-espace.

%%%%%%%%%%%%%%%%%%%%%%%%%%%%%%%%%%%%%
\section{Notations et rappels}\label{Rappels}

%%%%%%%%%%%%%%%%%%%%%%%%%
\subsection{Arbre de Bruhat-Tits, décompositions et projections de Cartan}\label{Cartan}

Soit $\kkk$ un corps local ultramétrique, c'est-à-dire une extension finie de~$\Q_p$ ou le corps $\F_q((t))$ des séries de Laurent formelles à coefficients dans un corps fini $\F_q$.
Fixons une valuation (additive)~$\omega$ sur~$\kkk$ à valeurs dans $\Z$.

Nous notons les $\kkk$-groupes algébriques par des lettres majuscules grasses (par exemple~$\mathbf{G}$) et leurs $\kkk$-points par la même lettre majuscule non grasse (par exemple~$G$).

%%%%%%%%
\subsubsection{Centralisateurs de tores et sous-groupes compacts ouverts de~$G$}

Soit $\mathbf{G}$ un $\kkk$-groupe algébrique semi-simple connexe de $\kkk$-rang~un.
Fixons un $\kkk$-tore $\kkk$-déployé maximal~$\mathbf{A}$ de~$\mathbf{G}$ et notons $\mathbf{Z}$ (resp.~$\mathbf{N}$) son centralisateur (resp. son normalisateur) dans~$\mathbf{G}$.
Soit $\pi$ un générateur du groupe des \linebreak $\kkk$-caractères de~$\mathbf{A}$ et soit $\alpha\in\N\pi$ une racine restreinte de~$\mathbf{A}$ dans~$\mathbf{G}$, c'est-à-dire un poids non trivial de $\mathbf{A}$ dans la représentation adjointe de~$\mathbf{G}$.
Quitte à remplacer~$\alpha$ par~$\alpha/2$, on peut supposer que~$\alpha$ est indivisible, c'est-à-dire que~$\alpha/2$ n'est pas une racine.
Pour tout $\kkk$-caractère~$\chi$ de~$\mathbf{Z}$, la restriction de~$\chi$ à~$\mathbf{A}$ est de la forme~$n_{\chi}\pi$, où~$n_{\chi}\in\Z$.
Pour tout~$z\in Z$ on note~$\nu(z)$ le réel défini par $n_{\chi}\nu(z) = -\omega(\chi(z))$ pour tout~$\chi$, et l'on pose
$$Z^+ = \{ z\in Z,\ \nu(z)\geq 0\} .$$
Le groupe~$Z$ agit sur~$\R$ par translation selon~$\nu$, et cette action se prolonge (de manière unique à une translation près) en une action affine de~$N$ sur~$\R$.
Notons~$\mathbf{U}_{\alpha}$ (resp.~$\mathbf{U}_{-\alpha}$) le $\kkk$-sous-groupe unipotent connexe de~$\mathbf{G}$ normalisé par~$\mathbf{Z}$ d'algèbre de Lie $\g_{\alpha}\oplus\g_{2\alpha}$ (resp.~$\g_{-\alpha}\oplus\g_{-2\alpha}$), où $\g_{i\alpha}$ désigne le sous-espace vectoriel de l'algèbre de Lie de~$\mathbf{G}$ formé des éléments~$X$ tels que $\Ad(a)(X)=\alpha(a)^iX$ pour tout~$a\in A$ (\cite{bor}, prop.~21.9).
Pour tout $u\in\nolinebreak U_{\alpha}\smallsetminus\nolinebreak\{ 1\} $, l'ensemble $N\cap U_{-\alpha}\,u\,U_{-\alpha}$ possède un unique élément, qui agit sur~$\R$ par la symétrie de centre~$x_u$ pour un certain~$x_u\in\R$.
Pour tout~$x\in\R$, posons
$$U_{\alpha,x} = \{ u\in U_{\alpha},\ u=1\ \mathrm{ou}\ x_u\leq x\} \ ;$$
d'après~\cite{bt2}, c'est un sous-groupe de~$U_{\alpha}$.
On définit de même un sous-groupe~$U_{-\alpha,x}$ de~$U_{-\alpha}$ pour tout~$x\in\R$.
Posons $N_x = \{ n\in N,\ n\cdot x=x\} $ et notons~$K_x$ le sous-groupe de~$G$ engendré par~$N_x$, $U_{\alpha,x}$ et~$U_{-\alpha,x}$.
Le groupe~$K_x$ est compact et ouvert dans~$G$, et maximal pour ces propriétés.

%%%%%%%%
\subsubsection{Arbre de Bruhat-Tits de~$G$}\label{Arbre de Bruhat-Tits}

On définit une relation d'équivalence sur~$G\times\R$ en posant $(g,x)\sim (g',x')$ s'il existe un élément~$n\in N$ tel que~$x'=n\cdot x$ et~$g^{-1}g'n\in K_x$.
Notons~$X$ l'ensemble des classes d'équivalence de~$G\times\R$ pour cette relation.
D'après~\cite{bt1} l'ensemble~$X$, muni de la topologie quotient induite par la topologie discrète de~$G$ et la topologie usuelle de~$\R$, est un arbre simplicial bipartite de valence~$\geq 3$, appelé \textit{arbre de Bruhat-Tits} de~$G$.
Il ne dépend pas du choix du $\kkk$-tore $\kkk$-déployé maximal~$\mathbf{A}$.
Le groupe~$G$ agit sur~$X$ par
$$g\cdot\overline{(g',x)} = \overline{(gg',x)}\,,$$
où~$\overline{(g,x)}$ désigne l'image de~$(g,x)\in G\times\R$ dans~$X$.
C'est une action propre, cocompacte, par automorphismes d'arbre simplicial.
Par construction, pour tout~$g\in\nolinebreak G$ et tout~$x\in\R$, le fixateur de~$\overline{(g,x)}$ dans~$G$ est le sous-groupe compact ouvert~$gK_xg^{-1}$.
D'après~\cite{bt1}, si l'on pose~$\mathcal{A}_0=\overline{\{ 1\} \times\R}$, alors les droites réelles plongées dans~$X$ sont exactement les ensembles $g\cdot\mathcal{A}_0$, où~$g\in G$.
On munit~$X$ de la distance~$d$ pour laquelle le plongement naturel de~$\{ g\} \times\R$ dans~$X$ est une isométrie pour tout~$g\in G$.
Le groupe~$G$ agit par isométries pour cette distance, toutes les arêtes de~$X$ ont la même longueur et les droites~$g\cdot\mathcal{A}_0$ sont les géodésiques de~$X$.
Nous renvoyons le lecteur aux textes fondateurs \cite{bt1} et~\cite{bt2} pour plus de détails, ainsi qu'à~\cite{rou} pour des explications plus élémentaires et à~\cite{ser}, \S~II.1, pour le cas de~$G=\nolinebreak\SL_2(\kkk)$.

%%%%%%%%
\subsubsection{Décompositions et projections de Cartan}

Notons~$K=K_{x_0}$ le fixateur dans~$G$ du sommet $x_0=\overline{(1,0)}$ de~$\mathcal{A}_0$.
D'après~\cite{bt1}, le groupe~$G$ agit transitivement sur l'ensemble des couples~$(\mathcal{A},\mathcal{C})$, où~$\mathcal{A}$ est une droite géodésique de~$X$ et~$\mathcal{C}$ une arête contenue dans~$\mathcal{A}$.
On en déduit la \textit{décomposition de Cartan} $G=KZ^+K$ : tout élément~$g\in G$ s'écrit $g=k_1zk_2$, où~$k_1,k_2\in K$ et~$z\in Z^+$.
Si~$g=k_1zk_2=k'_1z'k'_2$, où $k'_1,k'_2\in K$ et~$z'\in Z^+$, et si~$g\notin K$, alors $k_1^{-1}k'_1\in K\cap Z$ et~$k'_2k_2^{-1}\in K\cap Z$.
Posons $\mu(g)=\nu(z)\in\R^+$ ; on a
\begin{equation}\label{mu distance}
\mu(g) = d(x_0,g\cdot x_0),
\end{equation}
donc~$\mu(g)$ ne dépend pas de l'écriture de~$g$ comme produit d'éléments de~$K$, $Z^+$ et~$K$.
L'application $\mu :\nolinebreak G\rightarrow\nolinebreak\R^+$ ainsi obtenue est continue, propre, bi-$K$-invariante ; son image est l'intersection de~$\R^+$ avec un réseau de~$\R$.
On dit que~$\mu$ est la \textit{projection de Cartan} associée à la décomposition de Cartan $G=KZ^+K$.
Lorsque~$\mathbf{G}$ est déployé sur~$\kkk$, on a~$\nu(Z)=\nu(A)$ et~$G=KA^+K$, où~$A^+=A\cap Z^+$.
Ce n'est pas le cas en général, contrairement à la situation des groupes réels.

%%%%%%%%
\subsubsection{Sous-additivité des projections de Cartan}

D'après~(\ref{mu distance}) on a
\begin{equation}\label{inegalite triangulaire pour mu}
\mu(gg')\, \leq\, \mu(g) + \mu(g')
\end{equation}
pour tous~$g,g'\in G$ et
\begin{equation}\label{mu de l'inverse}
\mu(g^{-1})\, =\, \mu(g)
\end{equation}
puisque~$G$ agit sur~$X$ par isométries.
Rappelons que les $\kkk$-tores $\kkk$-déployés maximaux de~$\mathbf{G}$ sont tous conjugués sur~$\kkk$ (\cite{bot}, th.~4.21).
D'après~(\ref{inegalite triangulaire pour mu}), si~$\mathbf{A}'=g\mathbf{A}g^{-1}$ (où~$g\in G$) est un autre $\kkk$-tore $\kkk$-déployé maximal de $\mathbf{G}$ et si $\mu' : G\rightarrow\R^+$ est une projection de Cartan associée à~$\mathbf{A}'$, il existe une constante $C>0$ telle que
$$|\mu(g)-\mu'(g)|\, \leq\, C$$
pour tout~$g\in G$.
Ainsi, à une constante additive près, l'application~$\mu$ ne dépend pas du choix de~$\mathbf{A}$.

%%%%%%%%%%%%%%%%%%%%%%%%%%%%%%%%%%
\subsection{Isométries d'un arbre réel simplicial}\label{Action d'un groupe sur un arbre}

Nous considérons ici des arbres simpliciaux au sens de~\cite{ser}, déf.~6.
Rappelons qu'un \textit{arbre réel simplicial}, ou $\R$\textit{-arbre simplicial}, est un arbre simplicial muni d'une distance pour laquelle l'image de tout chemin injectif est isométrique à un segment réel.
Toutes les arêtes n'ont pas nécessairement la même longueur pour cette distance.

Fixons un arbre réel simplicial~$X$ de valence~$\geq 2$, notons~$d$ la distance sur~$X$ et~$\Isom(X)$ le groupe des isométries bijectives de~$X$ envoyant sommet sur sommet et arête sur arête.
Cette dernière condition est toujours vérifiée si~$X$ est de valence~$\geq 3$.
On munit~$\Isom(X)$ de la topologie pour laquelle les fixateurs (point par point) des parties compactes de~$X$ forment une base de voisinages compacts ouverts de l'identité.
Le groupe~$\Isom(X)$ est localement compact pour cette topologie.

%%%%%%%%
\subsubsection{Isométries hyperboliques et elliptiques}

Un élément~$g\in\Isom(X)$ est sans point fixe dans~$X$ si et seulement s'il existe une droite géodésique~$\mathcal{A}_g$, stable par~$g$, sur laquelle $g$ agit par une translation non triviale (\cite{ti1}, prop.~3.2, ou \cite{ser}, prop.~25) ; on note alors~$\lambda(g)>0$ la longueur de cette translation et l'on dit que~$g$ est \textit{hyperbolique}, d'\textit{axe de translation}~$\mathcal{A}_g$.
La figure~1 illustre l'action de~$g$ sur~$X$ dans ce cas.
Comme~$X$ est un arbre, pour tout~$x\in X$, le projeté~$\pr_g(x)$ de~$x$ sur~$\mathcal{A}_g$ est bien défini : c'est l'unique point de~$\mathcal{A}_g$ dont la distance à~$x$ est minimale.
On a
$$\pr_g(g\cdot x) = g\cdot\pr_g(x) \quad\quad\mathrm{et}\quad\quad d\big(g\cdot x,\mathcal{A}_g\big) = d\big(x,\mathcal{A}_g\big).$$
Les points $x$, $\pr_g(x)$, $g\cdot\pr_g(x)$ et~$g\cdot x$ sont alignés dans cet ordre sur une même droite géodésique de~$X$, d'où
\begin{eqnarray}\label{mu en fonction de lambda}
d(x,g\cdot x) & = & d\big(x,\pr_g(x)\big) + d\big(\pr_g(x),g\cdot\pr_g(x)\big) + d\big(g\cdot\pr_g(x),g\cdot x\big)\nonumber\\
& = & \lambda(g) + 2\,d(x,\mathcal{A}_g).
\end{eqnarray}
En particulier on a
$$\lambda(g)\ =\ \min_{x\in X}\, d(x,g\cdot x) > 0.$$

\begin{figure}\label{Hyperbolique}
\begin{center}\setlength{\unitlength}{1mm}
\begin{picture}(60,20)(-27,-4)
\put(-30,0){\line(1,0){60}}

\put(-15,0){\circle*{1}}
\put(-15,0){\line(-1,2){5}}
\put(-20,10){\circle*{1}}

\put(15,0){\circle*{1}}
\put(15,0){\line(1,2){5}}
\put(20,10){\circle*{1}}

\put(-22,12){$x$}
\put(-18,-4){$\mathrm{pr}_g(x)$}
\put(10,-4){$g\cdot\mathrm{pr}_g(x)$}
\put(18,12){$g\cdot x$}
\put(32,-1){$\mathcal{A}_g$}

\multiput(-5,2)(1,0){10}{\line(1,0){0.3}}
\put(4,2){\vector(1,0){1}}
\end{picture}
\end{center}
\caption{\textit{Action d'un élément hyperbolique}}
\end{figure}
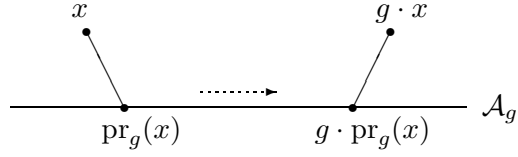

\noindent
Un élément~$g\in\Isom(X)$ qui admet un point fixe dans~$X$ est dit \textit{elliptique} ; l'ensemble~$X_g$ de ses points fixes est alors un sous-arbre de~$X$.
La figure~2 illustre l'action de~$g$ sur~$X$ dans ce cas.
Si pour tout~$x\in X$ on note~$\pr_g(x)$ le projeté de~$x$ sur~$X_g$, alors $d(\pr_g(x),x)=d(\pr_g(x),g\cdot x)$.
De plus, l'intersection des segments géodésiques~$[\pr_g(x),x]$ et~$[\pr_g(x),g\cdot x]$ est réduite à~$\{ \pr_g(x)\} $, d'où
\begin{equation}\label{distance de x a gx pour g non hyperbolique}
d(x,g\cdot x) = d(x,\pr_g(x)) + d(\pr_g(x),g\cdot x) = 2\,d(x,X_g).
\end{equation}
Par analogie avec le cas hyperbolique, on pose
$$\lambda(g)\ =\ \min_{x\in X}\, d(x,g\cdot x) = 0.$$

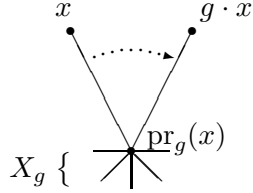
\begin{figure}\label{Elliptique}
\begin{center}\setlength{\unitlength}{1mm}
\begin{picture}(60,25)(-30,-4)
\put(0,0){\circle*{1}}

\put(0,0){\line(-1,2){8}}
\put(-8,16){\circle*{1}}

\put(0,0){\line(1,2){8}}
\put(8,16){\circle*{1}}

\put(-5,0){\line(1,0){10}}
\put(0,0){\line(-1,-1){4}}
\put(0,0){\line(0,-1){5}}
\put(0,0){\line(1,-1){4}}

\put(-10,18){$x$}
\put(2,1){$\mathrm{pr}_g(x)$}
\put(9,18){$g\cdot x$}
\put(-10,-3.3){$\big\{ $}
\put(-16,-3.5){$X_g$}

\put(-4.86,12.5){\circle*{0.5}}
\put(-3.68,12.9){\circle*{0.5}}
\put(-2.47,13.19){\circle*{0.5}}
\put(-1.24,13.36){\circle*{0.5}}
\put(0,13.42){\circle*{0.5}}
\put(1.24,13.36){\circle*{0.5}}
\put(2.47,13.19){\circle*{0.5}}
\put(3.68,12.9){\circle*{0.5}}
\put(4.86,12.5){\circle*{0.5}}
\put(4.66,12.55){\vector(4,-1){1}}
\end{picture}
\end{center}
\caption{\textit{Action d'un élément elliptique}}
\end{figure}

%%%%%%%%
\subsubsection{Bord de l'arbre~$X$}\label{Bord}

Comme tout espace hyperbolique au sens de Gromov, l'arbre~$X$ admet un \textit{bord}~$\partial X$, défini comme l'ensemble des classes d'équivalence de demi-droites géodésiques de~$X$ pour la relation ``être à distance de Hausdorff finie'' ou, de manière équivalente ici, ``être égales en dehors d'un compact''.
Si~$\xi\in\partial X$ désigne la classe d'une demi-droite géodésique~$D$, nous dirons que~$\xi$ est l'\textit{extrémité à l'infini} de~$D$.
Fixons un point~$x_0\in X$ et un réel~$q>1$.
Comme $X$ est un arbre, pour tout $\xi\in\partial X$ il existe une unique demi-droite géodésique~$D_{x_0}(\xi)$ d'extrémités $x_0\in X$ et~$\xi\in\partial X$.
Pour tous~$\xi,\xi'\in\partial X$ on pose $\dd(\xi,\xi') = q^{-r}$, où
$$r = \sup\big\{ d(x_0,x),\ x\in D_{x_0}(\xi)\cap D_{x_0}(\xi')\big\} \in [0,+\infty].$$
Ceci définit une distance~$\dd$ sur~$\partial X$.
L'action naturelle de~$\Isom(X)$ sur~$\partial X$ est continue pour cette distance.
Le fixateur de~$x_0$ agit sur~$\partial X$ par isométries.
Un élément~$g\in\Isom(X)$ est hyperbolique si et seulement s'il agit sur~$\partial X$ avec exactement deux points fixes, l'un attractif, noté~$\xi_g^+$, et l'autre répulsif, noté~$\xi_g^-$ ; ces points fixes sont les deux extrémités de l'axe de translation~$\mathcal{A}_g$.
Si~$g\in\Isom(X)$ est elliptique et fixe~$\xi\in\partial X$, alors $g$ fixe (point par point) toute une demi-droite de~$X$ d'extrémité~$\xi$.

%%%%%%%%
\subsubsection{Groupes discrets sans torsion d'isométries de~$X$}\label{groupes discrets d'isometries}

Un sous-groupe discret $\Gamma_0$ de~$\Isom(X)$ agit librement sur~$X$ si et seulement s'il est sans torsion.
Dans ce cas $\Gamma_0$ est un groupe libre (\cite{ser}, th.~4) et tous ses éléments non triviaux sont hyperboliques.
Plus précisément, pour tout domaine fondamental connexe~$\D$ de~$X$ pour l'action de~$\Gamma_0$, si l'on pose
$$\mathcal{F} = \big\{ \gamma\in\Gamma_0\smallsetminus\{ 1\} ,\quad \gamma\cdot\D\cap\D\neq\emptyset\big\} $$
et si~$\mathcal{F}'$ désigne une partie de~$\mathcal{F}$ telle que~$\mathcal{F}$ soit l'union disjointe de~$\mathcal{F}'$ et de~${\mathcal{F}'}^{-1}$, alors~$\mathcal{F}'$ est une partie libre génératrice de~$\Gamma_0$ (\cite{ser}, th.~$4'$).
Si~$\Gamma_0$ est de type fini, c'est un groupe de Schottky au sens de~\cite{lub}, déf.~1.4.
L'union~$X_{\Gamma}$ des axes de translations~$\mathcal{A}_{\gamma}$, où~$\gamma\in\Gamma_0\smallsetminus\{ 1\} $, est un sous-arbre de~$X$.
On dit que l'action de~$\Gamma_0$ sur~$X$ est \textit{minimale} si~$X_{\Gamma}=X$.
Lorsque $\Gamma_0$ est de type fini, son action sur~$X$ est minimale si et seulement si le graphe~$\Gamma_0\backslash X$ est fini (\cite{bas}, prop.~7.9, et \cite{bl}, th.~9.7).
Rappelons que l'\textit{ensemble limite} de~$\Gamma_0$ dans~$\partial X$ est par définition l'adhérence dans~$\partial X$ de l'ensemble des points fixes~$\xi_{\gamma}^+$, où $\gamma\in\Gamma_0\smallsetminus\{ 1\} $.
Si l'action de~$\Gamma_0$ sur~$X$ n'est pas minimale, l'ensemble limite de~$\Gamma_0$ dans~$\partial X$ est un fermé strict de~$\partial X$.
Ceci nous sera utile dans la démonstration du lemme~\ref{Proprietes des elements hyperboliques de G}.

%%%%%%%%%%%%%%%%%%%%%%%%%%%%%%%%%%%%%%%%%%%%%%%%%%%
\section{Une condition nécessaire et suffisante de cocompacité}\label{Cocompacite}

Cette partie est consacrée à la démonstration du théorème~\ref{theoreme quotients compacts}.
Dans toute la partie, nous notons~$\kkk$ un corps local ultramétrique, $G$ l'ensemble des \linebreak $\kkk$-points d'un $\kkk$-groupe algébrique semi-simple connexe de $\kkk$-rang~un et $\Delta_G$ la diagonale de~$G\times G$.
Nous fixons une décomposition de Cartan $G=KZ^+K$ et notons $\mu : G\rightarrow\R^+$ la projection de Cartan associée.

D'après \cite{kas}, th.~1.3, les sous-groupes discrets sans torsion de~$G\times G$ agissant proprement sur~$(G\times G)/\Delta_G$ sont, à la permutation près des deux facteurs de~$G\times G$, les graphes de la forme
$$\Gamma = \big\{ (\gamma,\rho(\gamma)),\ \gamma\in\Gamma_0\big\} ,$$
où $\Gamma_0$ est un sous-groupe discret sans torsion de~$G$ et $\rho : \Gamma_0\rightarrow G$ un morphisme de groupes admissible (au sens du théorème~\ref{theoreme quotients compacts}).
Pour démontrer le théorème~\ref{theoreme quotients compacts}, il suffit donc d'établir le résultat suivant.

\begin{theo}\label{equivalence cocompacite}
Soient $\kkk$ un corps local ultramétrique, $G$ l'ensemble des \linebreak $\kkk$-points d'un $\kkk$-groupe algébrique semi-simple connexe de $\kkk$-rang~un et $\Delta_G$ la diagonale de~$G\times G$.
Soient $\Gamma_0$ un sous-groupe discret de type fini sans torsion de~$G$ et $\rho : \Gamma_0\rightarrow G$ un morphisme de groupes admissible.
Notons
$$\Gamma = \big\{ (\gamma,\rho(\gamma)),\ \gamma\in\Gamma_0\big\} $$
le graphe de~$\rho$.
Le quotient~$\Gamma\backslash (G\times G)/\Delta_G$ est compact si et seulement si le quotient~$\Gamma_0\backslash G$ l'est.
\end{theo}

Nous utilisons pour cela l'existence d'un isomorphisme
\begin{eqnarray*}
(\mathbf{G}\times\mathbf{G})/\mathbf{\Delta_G} & \longrightarrow & \ \mathbf{G}\\
(g,h)\,\mathbf{\Delta_G}\ \, & \longmapsto & gh^{-1}.
\end{eqnarray*}
de $(\mathbf{G}\times\mathbf{G})$-variétés algébriques sur~$\kkk$, où $\mathbf{G}\times\mathbf{G}$ agit sur~$(\mathbf{G}\times\mathbf{G})/\mathbf{\Delta_G}$ par translation à gauche et sur~$\mathbf{G}$ par
$$(g_1,g_2)\cdot g\ =\ g_1\,g\,g_2^{-1}.$$
Cet isomorphisme induit un isomorphisme de $(G\times G)$-ensembles sur les \linebreak $\kkk$-points.

Notons $X$ l'arbre de Bruhat-Tits de~$G$, muni de la distance~$d$ définie au paragraphe~\ref{Arbre de Bruhat-Tits}, et~$\partial X$ son bord, muni de la distance~$\dd$ définie au paragraphe~\ref{Bord} en prenant par exemple pour~$q$ le cardinal du corps résiduel de~$\kkk$.
Pour traiter l'implication directe du théorème~\ref{equivalence cocompacite}, nous introduisons certains points~$\zeta_g^-$ du bord de l'arbre de Bruhat-Tits de~$G$, obtenus à partir de la décomposition de Cartan~$G=KZ^+K$.

%%%%%%%%%%%%%%%%%%%%%%%%%
\subsection{Points de~$\partial X$ associés à la décomposition de Cartan~$G=KZ^+K$}\label{Points y}

Pour tout élément hyperbolique~$g\in G$, notons~$\xi_g^+$ (resp.~$\xi_g^-$) son point fixe attractif (resp. répulsif) dans~$\partial X$.
Notons~$\xi_{Z^+}^-$ le point fixe répulsif commun à tous les éléments hyperboliques de~$Z^+$.

Soit $g\in G$ de décomposition de Cartan $g = k_1 z k_2$, où $k_1,k_2\in K$ et $z\in Z^+$.
Si~$g\notin K$, le point
\begin{equation}\label{definition de y}
\zeta_g^- = k_2^{-1}\cdot\xi_{Z^+}^-\in\partial X
\end{equation}
est bien défini, car si~$g=k'_1z'k'_2$ est une autre décomposition de Cartan de~$g$, où~$k'_1,k'_2\in K$ et~$z'\in Z^+$, alors~${k'}_2^{-1}\in k_2^{-1}Z$ (\textit{cf.} paragraphe~\ref{Cartan}).
Le lemme suivant montre que lorsque~$g$ est hyperbolique et $\mu(g)$ grand, le point~$\zeta_g^-$ est proche du point fixe répulsif~$\xi_g^-$.

\begin{lem}\label{y proche de x}
Pour tout élément hyperbolique~$g\in G$, on a
$$\dd(\xi_g^-,\zeta_g^-)\, \leq\, q^{-\mu(g)/2}.$$
\end{lem}

\begin{dem}
Soit $g\in G$ un élément hyperbolique ; écrivons~$g=k_1zk_2$, où~$k_1,k_2\in K$ et~$z\in Z^+$.
Les points $x_0$, $\pr_{g^{-1}}(x_0)$, $g^{-1}\cdot\pr_{g^{-1}}(x_0)$ et $g^{-1}\cdot x_0$ sont alignés dans cet ordre, donc $g^{-1}\cdot\pr_{g^{-1}}(x_0) \in D_{x_0}(\xi_g^-)$.
D'autre part, comme $K$ fixe~$x_0$, on a
$$g^{-1}\cdot x_0 = k_2^{-1}z^{-1}\cdot x_0\, \in\, k_2^{-1}\cdot D_{x_0}(\xi_{Z^+}^-) = D_{x_0}(\zeta_g^-),$$
donc $g^{-1}\cdot\pr_{g^{-1}}(x_0)\in D_{x_0}(\zeta_g^-)$ car~$X$ est un arbre.
On en déduit
$$\dd(\xi_g^-,\zeta_g^-) \leq q^{-d(x_0,g^{-1}\cdot\pr_{g^{-1}}(x_0))}.$$
Or, d'après~(\ref{mu distance}), (\ref{mu de l'inverse}) et~(\ref{mu en fonction de lambda}) on a
$$d\big(x_0,g^{-1}\cdot\pr_{g^{-1}}(x_0)\big)\ \leq\ \frac{\mu(g)}{2},$$
d'où le résultat.
\end{dem}

\bigskip

L'intérêt d'introduire les points~$\zeta_g^-$ est de contrôler la projection de Cartan de certains éléments de $G$, comme dans le lemme suivant.

\medskip

\begin{lem}\label{Interet des points y}
Pour tout élément hyperbolique~$g\in G$ et tout réel~$\varepsilon>0$, il existe des constantes~$C,N\geq 0$ telles que pour tout $\gamma\in G\smallsetminus K$ vérifiant $\dd(\xi_g^+,\zeta_{\gamma}^-)\geq\varepsilon$ on ait
$$\mu(\gamma g^n) \geq \mu(\gamma) + \mu(g^n) - C$$
pour tout entier $n>N$.
\end{lem}

\begin{dem}
Soient~$\varepsilon>0$ un réel et~$g\in G$ un élément hyperbolique.
Posons
$$s = -\log_q\varepsilon$$
et montrons que $C=2s$ et $N=s/\lambda(g)$ conviennent.
La figure~3 illustre notre raisonnement.
Pour tout~$n\geq\nolinebreak 1$, les points $x_0$, $\pr_g(x_0)$ et~$g^n\cdot\pr_g(x_0)$ sont alignés dans cet ordre sur la demi-droite géodésique~$D_{x_0}(\xi_g^+)$, d'où
\begin{equation}\label{x et x_0}
d\big(x_0,g^n\cdot\pr_g(x_0)\big)\ \geq\ d\big(\pr_g(x_0),g^n\cdot\pr_g(x_0)\big)\ =\ n\lambda(g).
\end{equation}
Soit~$\gamma\in G\smallsetminus K$ tel que $\dd(\xi_g^+,\zeta_{\gamma}^-)\geq\varepsilon$.
\'Ecrivons $\gamma = k_1 z k_2$, où~$k_1,k_2\in K$ et~$z\in Z^+$.
L'intersection~$D_{x_0}(\xi_g^+)\cap D_{x_0}(\zeta_{\gamma}^-)$ est un segment géodésique dont l'une des extrémités est~$x_0$ ; notons~$x\in X$ son autre extrémité.
On a $d(x_0,x)\leq s$ par hypothèse.
Soit $n > s/\lambda(g)$ un entier.
D'après~(\ref{x et x_0}) on~a
$$d(x_0,g^n\cdot\pr_g(x_0))\, >\, d(x_0,x),$$
donc les points $x_0$, $x$ et~$g^n\cdot x_0$ sont alignés dans cet ordre et par~(\ref{mu distance}) on a
\begin{equation}\label{calcul distance 1}
d(x,g^n\cdot x_0)\ =\ \mu(g^n) - d(x_0,x).
\end{equation}

\begin{figure}\label{Figure lemme}
\begin{center}\setlength{\unitlength}{1mm}
\begin{picture}(95,56)(0,-12)
\put(0,0){\line(1,0){85}}

\put(0,0){\circle*{1}}
\put(-1,-4){$x_0$}

\put(20,0){\circle*{1}}
\put(16,-4){$\mathrm{pr}_g(x_0)$}

\put(35,0){\circle*{1}}
\put(34,-4){$x$}

\put(60,0){\circle*{1}}
\put(56,-4){$g^n\cdot\mathrm{pr}_g(x_0)$}

\put(20,0){\line(-1,2){3}}
\put(16,8){\vector(-1,2){2}}
\put(11,14){$\xi_g^-$}

\put(35,0){\line(1,2){15}}
\put(51,32){\vector(1,2){2.3}}
\put(55,38){$\zeta_{\gamma}^-=k_2^{-1}\cdot\xi_{Z^+}^-$}

\put(45,20){\circle*{1}}
\put(47,19){$\gamma^{-1}\cdot x_0$}

\put(60,0){\line(1,1){14.14}}
\put(74.14,14.14){\circle*{1}}
\put(76,13){$g^n\cdot x_0$}

\put(87,0){\vector(1,0){5}}
\put(95,-1){$\xi_g^+$}

\multiput(0,-8)(1,0){35}{\line(1,0){0.3}}
\put(1,-8){\vector(-1,0){1}}
\put(34,-8){\vector(1,0){1}}
\put(14,-12){$\leq s$}
\end{picture}
\end{center}
\caption{\textit{Illustration du lemme~\ref{Interet des points y}}}
\end{figure}
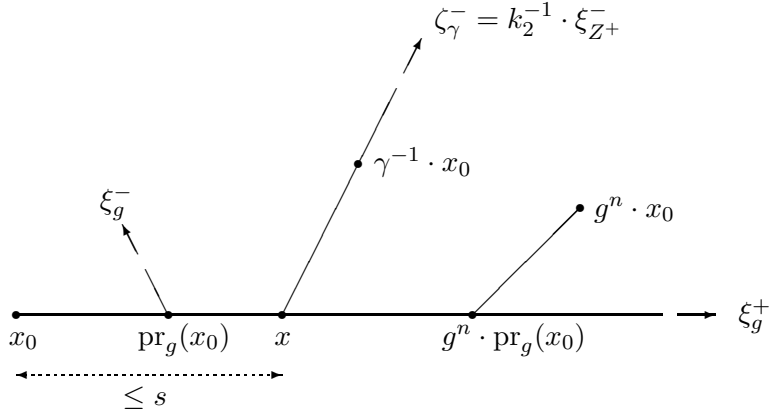

\noindent
Comme $K$ fixe~$x_0$, on a $\gamma^{-1}\cdot x_0=k_2^{-1}z^{-1}\cdot x_0$, et ce point appartient à $k_2^{-1}\cdot D_{x_0}(\xi_{Z^+}^-)=D_{x_0}(\zeta_{\gamma}^-)$.
Le point~$x$ appartient lui aussi à~$D_{x_0}(\zeta_{\gamma}^-)$, donc
$$d(\gamma^{-1}\cdot x_0,x) = \big|d(x_0,\gamma^{-1}\cdot x_0) - d(x_0,x)\big|.$$
En utilisant (\ref{mu distance}) et~(\ref{mu de l'inverse}), on obtient
\begin{equation}\label{calcul distance 2}
d(\gamma^{-1}\cdot x_0,x)\ \geq\ \mu(\gamma) - d(x_0,x).
\end{equation}
Or, $x$ est le projeté de~$g^n\cdot x_0$ sur~$D_{x_0}(\zeta_{\gamma}^-)$, donc les points $\gamma^{-1}\cdot x_0$, $x$ et~$g^n\cdot x_0$ sont alignés dans cet ordre, d'où
\begin{equation}\label{calcul distance 3}
d(\gamma^{-1}\cdot x_0,g^n\cdot x_0)\ =\ d(\gamma^{-1}\cdot x_0,x) + d(x,g^n\cdot x_0).
\end{equation}
En utilisant (\ref{mu distance}), (\ref{calcul distance 1}), (\ref{calcul distance 2}), (\ref{calcul distance 3}) et le fait que~$G$ agit sur~$X$ par isométries, on trouve
\begin{eqnarray*}
\mu(\gamma g^n) & = & d(x_0,\gamma g^n\cdot x_0)\\
& = & d(\gamma^{-1}\cdot x_0,g^n\cdot x_0)\\
& \geq & \mu(\gamma) + \mu(g^n) - 2\,d(x_0,x),
\end{eqnarray*}
d'où le résultat puisque $d(x_0,x)\leq s$.
\end{dem}

%%%%%%%%%%%%%%%%%%%%%%%%%%%%%%%%%%
\subsection{Démonstration de l'implication directe du théorème \ref{equivalence cocompacite}}\label{Sens direct}

Pour obtenir l'implication directe du théorème~\ref{equivalence cocompacite}, nous raisonnons par contraposition : il suffit de démontrer la proposition suivante.

\begin{prop}\label{sens direct par contraposition}
Soient $\Gamma_0$ un sous-groupe discret de type fini sans torsion de~$G$ et $\rho : \Gamma_0\rightarrow G$ un morphisme de groupes admissible.
Si $\Gamma_0$ n'est pas cocompact dans~$G$, alors pour tout~$R>0$ il existe un élément~$g_R\in G$ tel que $\mu(\gamma\,g_R\,\rho(\gamma)^{-1}) \geq R$ pour tout~$\gamma\in\Gamma_0$.
\end{prop}

Pour démontrer la proposition~\ref{sens direct par contraposition}, nous établissons l'existence d'un élément hyperbolique~$g\in G$ tel que pour presque tout~$\gamma\in\Gamma_0\smallsetminus\{ 1\} $, le point $\zeta_{\gamma}^-\in\partial X$ donné par~(\ref{definition de y}) soit suffisamment éloigné du point fixe attractif~$\xi_g^+$ (lemme~\ref{Proprietes des elements hyperboliques de G}).
Nous utilisons ensuite le lemme~\ref{Interet des points y}.

\medskip

\begin{lem}\label{Proprietes des elements hyperboliques de G}
Soit $\Gamma_0$ un sous-groupe discret de type fini sans torsion de~$G$, non cocompact dans~$G$.
Il existe un élément hyperbolique~$g\in G$ et un réel $\varepsilon>\nolinebreak 0$ tels que $\dd(\xi_g^+,\zeta_{\gamma}^-) \geq \varepsilon$ pour presque tout~$\gamma\in\nolinebreak\Gamma_0\smallsetminus\nolinebreak\{ 1\} $.
\end{lem}

\begin{dem}
Comme~$\Gamma_0$ est discret, sans torsion, de type fini et non cocompact dans~$G$, son ensemble limite dans~$\partial X$ est un fermé strict de~$\partial X$ (\textit{cf.} paragraphe~\ref{groupes discrets d'isometries}).
Notons~$\mathcal{U}$ le complémentaire dans~$\partial X$ de cet ensemble limite.
Soit $g\in G$ un élément hyperbolique tel que~$\xi_g^+\in\mathcal{U}$.
Il existe un réel~$\varepsilon>0$ tel que
$$\big\{ \xi\in\partial X,\ \dd\big(\xi_g^+,\xi\big)<2\varepsilon\big\} \ \subset\ \mathcal{U}.$$
Par définition on a alors $\dd(\xi_g^+,\xi_{\gamma}^-)\geq 2\varepsilon$ pour tout~$\gamma\in\Gamma_0\smallsetminus\{ 1\} $.
Le lemme~\ref{y proche de x} implique
$$\dd\big(\xi_g^+,\zeta_{\gamma}^-\big)\ \geq\ \dd\big(\xi_g^+,\xi_{\gamma}^-\big) - \dd\big(\xi_{\gamma}^-,\zeta_{\gamma}^-\big)\ \geq\ 2\varepsilon - q^{-\mu(\gamma)/2}$$
pour tout~$\gamma\in\Gamma_0\smallsetminus\{ 1\} $.
Or, l'application~$\mu$ est propre et le groupe~$\Gamma_0$ discret dans~$G$, donc $q^{-\mu(\gamma)/2}\leq\varepsilon$ pour presque tout~$\gamma\in\Gamma_0\smallsetminus\{ 1\} $.
\end{dem}

\bigskip\smallskip

Nous pouvons à présent démontrer la proposition~\ref{sens direct par contraposition}.

\medskip

\noindent
\textbf{Démonstration de la proposition~\ref{sens direct par contraposition}.}
Supposons $\Gamma_0$ non cocompact dans~$G$.
D'après les lemmes~\ref{Interet des points y} et~\ref{Proprietes des elements hyperboliques de G}, il existe un élément hyperbolique~$g\in G$ et des réels~$C,N\geq 0$ tels que pour tout~$\gamma\in\Gamma_0$ en dehors d'un certain ensemble fini~$F$, on ait
\begin{equation}\label{minoration de mu}
\mu(\gamma g^n)\ \geq\ \mu(\gamma)+\mu(g^n)-C
\end{equation}
pour tout entier~$n>N$.
Soit $R>0$.
Comme~$\rho$ est admissible, on a
\begin{equation}\label{majoration de mu de rho}
\mu(\rho(\gamma))\ \leq\ \mu(\gamma)-R-C
\end{equation}
pour tout~$\gamma\in\Gamma_0$ en dehors d'un certain ensemble fini~$F'$. Posons
$$R' = \max_{\gamma\in F\cup F'} \big(\mu(\gamma^{-1})+\mu(\rho(\gamma))\big).$$
Comme $g$ est hyperbolique, la suite $(\mu(g^n))_{n\in\N}$ tend vers l'infini avec~$n$ d'après~(\ref{mu en fonction de lambda}).
En particulier, il existe un entier~$n>N$ tel que~$\mu(g^n)\geq\nolinebreak R+\nolinebreak R'$.
Fixons un tel~$n>N$ et montrons que l'élément~$g_R=g^n$ convient.
D'après (\ref{inegalite triangulaire pour mu}), (\ref{minoration de mu}) et~(\ref{majoration de mu de rho}), pour tout~$\gamma\notin (F\cup F')$ on a
$$\mu\big(\gamma g^n \rho(\gamma)^{-1}\big) \geq \mu(\gamma g^n) - \mu(\rho(\gamma)) \geq \big(\mu(\gamma) + \mu(g^n) - C\big) - \big(\mu(\gamma) - R - C\big) \geq R$$
et pour tout $\gamma\in F\cup F'$ on a
$$\mu\big(\gamma g^n \rho(\gamma)^{-1}\big)\ \geq\ \mu(g^n) - \mu(\gamma^{-1}) - \mu(\rho(\gamma))\ \geq\ \mu(g^n) - R'\ \geq\ R.$$
Ceci achève la démonstration de la proposition~\ref{sens direct par contraposition}.
\hfill\qedsymbol

%%%%%%%%%%%%%%%%%%%%%%%%%%%%%%%%%
\subsection{Démonstration de l'implication réciproque du théorème \ref{equivalence cocompacite}}\label{Sens reciproque}

Pour terminer la démonstration du théorème~\ref{equivalence cocompacite} (et donc du théorème~\ref{theoreme quotients compacts}), il suffit d'établir le résultat suivant.

\begin{prop}\label{proposition sens reciproque}
Soient $\Gamma_0$ un réseau cocompact sans torsion de~$G$ et $\rho :\nolinebreak \Gamma_0\rightarrow G$ un morphisme de groupes admissible.
Il existe une partie compacte~$\mathcal{C}$ de~$G$ telle que
$$G = \big\{ \gamma g \rho(\gamma)^{-1},\quad \gamma\in\Gamma_0,\ g\in\mathcal{C}\big\} .$$
\end{prop}

\begin{dem}
Par cocompacité de~$\Gamma_0$ dans~$G$ et continuité de~$\mu$, il existe un réel~$R>0$ tel que
$$G = \Gamma_0\cdot\big\{ g\in G,\ \mu(g)\leq R\big\} .$$
Comme~$\rho$ est admissible, on a
$$\mu(\rho(\gamma))\ \leq\ \mu(\gamma)-2R-1$$
pour tout~$\gamma\in\Gamma_0$ en dehors d'un certain ensemble fini~$F$.
Posons $R' = \max_{\gamma\in F} \mu(\gamma)$ et montrons que le compact
$$\mathcal{C} = \big\{ g\in G,\ \mu(g) \leq R + R'\big\} $$
convient.
Pour cela, il suffit d'établir que pour tout élément~$g\in G$ vérifiant $\mu(g)>R+R'$ il existe~$\gamma\in\Gamma_0$ tel que~$\mu(\gamma g \rho(\gamma)^{-1})\leq\mu(g)-1$ ; on peut alors conclure par récurrence.

Soit~$g\in G$ tel que~$\mu(g)>R+R'$.
Par hypothèse, il existe~$\gamma\in\Gamma_0$ tel que~$\mu(\gamma g)\leq\nolinebreak R$.
Montrons que $\mu(\gamma g \rho(\gamma)^{-1})\leq\mu(g)-1$.
D'après~(\ref{inegalite triangulaire pour mu}) on a
$$\mu(\gamma^{-1})\ \geq\ \mu(g) - \mu(\gamma g)\ >\ R',$$
donc $\gamma^{-1}\notin F$, ce qui implique $\mu(\rho(\gamma)^{-1})\leq\mu(\gamma^{-1})-2R-1$.
En utilisant~(\ref{inegalite triangulaire pour mu}) et~(\ref{mu de l'inverse}), on trouve
$$\mu(\gamma^{-1})\ \leq\ \mu(g) + \mu(g^{-1}\gamma^{-1})\ \leq\ \mu(g) + R$$
et
$$\mu\big(\gamma g \rho(\gamma)^{-1}\big)\ \leq\ \mu(\gamma g) + \mu\big(\rho(\gamma)^{-1}\big)\ \leq\ \mu(\gamma^{-1}) - R - 1\ \leq\ \mu(g) - 1.$$
Ceci achève la démonstration de la proposition~\ref{proposition sens reciproque}.
\end{dem}

%%%%%%%%%%%%%%%%%%%%%%%%%%%%%%%%%
\subsection{Le cas de rang supérieur}\label{Rang superieur}

La description donnée par le théorème~\ref{theoreme quotients compacts} est spécifique au rang~un.
En effet, voici en caractéristique nulle un exemple de quotient compact $\Gamma\backslash (G\times G)/\Delta_G$ où $\mathbf{G}$ est de $\kkk$-rang~$\geq 2$ et où $\Gamma$ est le produit de deux sous-groupes infinis de~$G$.
Fixons un élément non carré $\beta\in\kkk\smallsetminus\kkk^2$.
Soient~$Q$ la forme quadratique sur~$\kkk^4$ donnée par
$$Q(x_1,y_1,x_2,y_2) = (x_1^2 - \beta\,y_1^2) - (x_2^2 - \beta\,y_2^2)$$
et $\mathbf{G}=\mathbf{SO}(Q)$ le groupe spécial orthogonal de~$Q$.
Notons~$Q_1$ la restriction de~$Q$ à~$\kkk^3\times\{ 0\} $ et soit~$\mathbf{H}_1=\mathbf{SO}(Q_1)$ le groupe spécial orthogonal de~$Q_1$, vu comme sous-groupe de~$\mathbf{G}$.
Choisissons une racine carrée~$\sqrt{\beta}$ de~$\beta$ dans une clôture algébrique de~$\kkk$ et notons $\sigma$ l'élément non trivial du groupe de Galois de l'extension quadratique~$\kkk(\sqrt{\beta})/\kkk$.
Soient~$h$ la forme hermitienne sur~$\kkk(\sqrt{\beta})^2$ donnée par
$$h(z_1,z_2) = z_1\,\sigma(z_1) - z_2\,\sigma(z_2)$$
et $\mathbf{H}_2=\mathbf{U}(h)$ le groupe unitaire de~$h$.
En identifiant~$\kkk(\sqrt{\beta})$ à~$\kkk^2$ par l'application qui à tout~$x+\sqrt{\beta}y$ associe~$(x,y)$, on voit~$\mathbf{H}_2$ comme un sous-groupe de~$\mathbf{G}$.
Par un analogue ultramétrique de \cite{ko1}, prop.~4.9, le groupe $H_1$ agit proprement et cocompactement sur~$G/H_2$, donc $H_1\times H_2$ agit proprement et cocompactement sur~$(G\times G)/\Delta_G$.
Si $\kkk$ est de caractéristique nulle, les groupes~$H_1$ et~$H_2$ admettent des réseaux cocompacts sans torsion~$\Gamma_1$ et~$\Gamma_2$ (\cite{bh}, th.~A, et \cite{sel}, lem.~8).
Le groupe $\Gamma=\Gamma_1\times\Gamma_2$ agit alors librement, proprement et cocompactement sur~$(G\times G)/\Delta_G$.

%%%%%%%%%%%%%%%%%%%%%%%%%%%%%%%%%%%%%%%%%%%%%%%%%%%%%%%%%%%%
\section{Longueurs de translation et constantes de Lipschitz}\label{Applications equivariantes entre arbres}

Reprenons les notations du paragraphe~\ref{Action d'un groupe sur un arbre}.
Le but de cette partie est de démontrer le résultat suivant.

\begin{prop}\label{max sur un ensemble fini}
Soient $X$ et~$X'$ deux arbres réels simpliciaux de valence~$\geq\nolinebreak 2$.
Soit~$\Gamma_0$ un sous-groupe discret sans torsion de~$\Isom(X)$ tel que le graphe quotient~$\Gamma_0\backslash X$ soit fini, et soit $\rho : \Gamma_0\rightarrow\nolinebreak\Isom(X')$ un morphisme de groupes.
\begin{enumerate}
	\item Il existe une application affine par morceaux $f : X\rightarrow~X'$ qui est $\rho$-équivariante, au sens où $f(\gamma\cdot x) = \rho(\gamma)\cdot f(x)$ pour tout~$\gamma\in\Gamma_0$ et tout~$x\in X$.
	Une telle application est lipschitzienne.
	\item La borne inférieure~$C\geq 0$ des constantes de Lipschitz de telles applications est atteinte.
	\item Fixons un point~$x_0\in X$ et soit~$F$ le sous-ensemble fini de~$\Gamma_0\smallsetminus\{ 1\} $ formé des éléments~$\gamma$ tels que $d(x_0,\gamma\cdot x_0)\leq 4N$, où $N\geq 1$ désigne le nombre de sommets du graphe fini~$\Gamma_0\backslash X$.
	On a
	$$\sup_{\gamma\in\Gamma_0\smallsetminus\{ 1\} } \frac{\lambda\big(\rho(\gamma)\big)}{\lambda(\gamma)}\ =\ \max_{\gamma\in F}\, \frac{\lambda\big(\rho(\gamma)\big)}{\lambda(\gamma)}\ =\ C.$$
\end{enumerate}
\end{prop}

\bigskip

Remarquons que l'hypothèse que $\Gamma_0$ est un sous-groupe discret sans torsion de~$\Isom(X)$ tel que le graphe~$\Gamma_0\backslash X$ soit fini est équivalente à l'hypothèse que $\Gamma_0$ est un groupe libre de type fini agissant sur~$X$ librement, minimalement, par isométries (\textit{cf.} paragraphe~\ref{groupes discrets d'isometries}).

Nous notons ici $d$ la distance sur~$X$ ou~$X'$ et disons qu'une application $f : X\rightarrow X'$ est \textit{affine par morceaux} si pour toute arête~$e$ de~$X$, il existe une constante~$C_e\geq 0$ telle que
\begin{equation}\label{affine}
d\big(f(x_1),f(x_2)\big) = C_e\,d(x_1,x_2)
\end{equation}
pour tous~$x_1,x_2\in e$.
Une telle application est continue.

\begin{rema}\label{remarque affine par morceaux}
Une application $f : X\rightarrow X'$ affine par morceaux est entièrement déterminée par les images par~$f$ des sommets de~$X$.
\end{rema}

Notons que pour tout~$\gamma\in\Gamma_0\smallsetminus\{ 1\} $, tout~$x\in\mathcal{A}_{\gamma}$ et toute application $\rho$-équivariante et $C'$-lipschitzienne $f : X\rightarrow X'$, on a
$$\lambda(\rho(\gamma))\, \leq\, d\big(f(x),\rho(\gamma)\cdot f(x)\big)\, =\, d\big(f(x),f(\gamma\cdot x)\big)\, \leq\, C' d(x,\gamma\cdot x)\, =\, C' \lambda(\gamma).$$
L'inégalité
\begin{equation}\label{inegalite sup-inf}
\sup_{\gamma\in\Gamma_0\smallsetminus\{ 1\} } \frac{\lambda\big(\rho(\gamma)\big)}{\lambda(\gamma)}\ \leq\ C.
\end{equation}
en résulte immédiatement.
Le point~(3) de la proposition~\ref{max sur un ensemble fini} affirme que cette inégalité est une égalité et que la borne supérieure dans le terme de gauche est atteinte sur un sous-ensemble fini~$F$ indépendant de~$\rho$.

Dans le cas particulier où~$\rho$ est injectif d'image discrète et cocompacte, ce résultat est équivalent à un résultat de T.~White obtenu en étudiant une certaine ``distance asymétrique'' sur l'outre-espace (\cite{fm}, prop.~3.15).
Nous nous inspirons de la preuve de White, qui nous a été aimablement communiquée par M.~Bestvina.

%%%%%%%%%%%%%%%%%%%%%%%%%%%%%%%%%
\subsection{Applications $\rho$-équivariantes affines par morceaux}

Le point~(1) de la proposition~\ref{max sur un ensemble fini} est facile.

\medskip

\noindent
\textbf{Démonstration du point}~(1) \textbf{de la proposition~\ref{max sur un ensemble fini}.}
Montrons tout d'abord qu'il existe une application $f : X\rightarrow X'$ $\rho$-équivariante et affine par morceaux.
Comme le graphe~$\Gamma_0\backslash X$ est fini, il existe un système fini~$S$ de représentants de l'ensemble des sommets de~$X$ \textit{modulo} l'action de~$\Gamma_0$.
Pour tout~$s\in S$, choisissons un sommet~$s'$ de~$X'$.
D'après la remarque~\ref{remarque affine par morceaux}, il existe une unique application affine par morceaux $f : X\rightarrow X'$ telle que $f(\gamma\cdot s)=\rho(\gamma)\cdot s'$ pour tout~$\gamma\in\Gamma_0$ et tout~$s\in S$.
Cette application est $\rho$-équivariante.

D'autre part, toute application $f : X\rightarrow X'$ $\rho$-équivariante et affine par morceaux est lipschitzienne.
En effet, comme le graphe~$\Gamma_0\backslash X$ est fini, il existe un système fini~$E$ de représentants de l'ensemble des arêtes de~$X$ \textit{modulo} l'action de~$\Gamma_0$.
L'application $f$ est lipschitzienne de constante~$\max_{e\in E} C_e$, où~$C_e$ est donnée par~(\ref{affine}).
\hfill\qedsymbol

%%%%%%%%%%%%%%%%%%%%%%%%%%%%%%%%%
\subsection{Une constante de Lipschitz minimale}\label{constante de Lipschitz minimale}

Dans le cas particulier où $\rho$ est injectif d'image discrète et cocompacte, toute application $\rho$-équivariante affine par morceaux $f : X\rightarrow X'$ induit une application affine par morceaux du graphe fini~$\Gamma_0\backslash X$ vers le graphe fini~$\rho(\Gamma_0)\backslash X'$ ; le point~(2) de la proposition~\ref{max sur un ensemble fini} est alors une conséquence immédiate du théorème d'Ascoli.

Le cas général est plus compliqué : un argument supplémentaire est nécessaire pour montrer que l'on peut se ramener à des applications à valeurs dans un compact de~$X'$ et appliquer le théorème d'Ascoli.
Nous allons distinguer plusieurs cas selon le nombre de points fixes du groupe~$\rho(\Gamma_0)$ dans~$\partial X'$.
Commençons par un lemme.

\begin{lem}\label{lemme axes de translation et points fixes}
Soient~$X'$ un arbre réel simplicial, $g_1,\ldots,g_n\in\Isom(X')$ des isométries de~$X'$ sans point fixe commun dans~$X'$ et~$R'\geq 0$ un réel.
Posons
$$\mathcal{C} = \big\{ x'\in X',\quad d(x',g_i\cdot x')\leq R'\quad \forall i\big\} .$$
\begin{itemize}
	\item Si les isométries~$g_i$ n'admettent pas de point fixe commun dans~$\partial X'$, alors~$\mathcal{C}$ est compact.
	\item Si les isométries~$g_i$ admettent exactement un point fixe commun~$\xi'$ dans~$\partial X'$, il existe une demi-droite géodésique~$D$ de~$X'$ d'extrémité~$\xi'$ telle que $d(x',D)\leq R'/2$ pour tout~$x'\in\mathcal{C}$.
	\item Si les isométries~$g_i$ admettent deux points fixes communs~$\xi'_1$ et~$\xi'_2$ \linebreak dans~$\partial X'$, on a $d(x',D)\leq R'/2$ pour tout~$x'\in\mathcal{C}$, où~$D$ désigne la droite géodésique de~$X'$ d'extrémités $\xi'_1$ et~$\xi'_2$.
\end{itemize}
\end{lem}

\begin{dem}
Pour tout~$i$, si~$g_i$ est hyperbolique, notons~$X'_i$ l'axe de translation de~$g_i$ dans~$X'$.
Si~$g_i$ est elliptique, notons~$X'_i$ l'ensemble des points fixes de~$g_i$ dans~$X'$ ; c'est un sous-arbre de~$X'$ et les points fixes de~$g_i$ dans~$\partial X'$ sont exactement les extrémités à l'infini des demi-droites géodésiques incluses dans~$X'_i$.
D'après~(\ref{mu en fonction de lambda}) et~(\ref{distance de x a gx pour g non hyperbolique}), on a $d(x',X'_i)\leq R'/2$ pour tout~$x'\in\mathcal{C}$ et tout~$1\leq i\leq n$.

Supposons que les~$g_i$ n'admettent pas de point fixe commun dans~$\partial X'$ et qu'ils soient tous elliptiques.
Alors l'intersection des~$X'_i$ est vide car les~$g_i$ n'admettent pas de point fixe commun dans~$X'$ par hypothèse.
Comme~$X'$ est un arbre, il existe deux entiers~$i_1$ et~$i_2$ tels que $X'_{i_1}\cap X'_{i_2}=\emptyset$.
Pour tout~$x'\in\mathcal{C}$, nous avons vu que $d(x',X'_i)\leq R'/2$ pour tout~$i$.
On en déduit facilement que $d(x',x'')\leq R'/2$ pour tout~$x'\in\mathcal{C}$, où $x''$ désigne le projeté de~$X'_{i_1}$ sur~$X'_{i_2}$.
Ainsi, le fermé~$\mathcal{C}$ est borné, donc compact.

Supposons que les~$g_i$ n'admettent pas de point fixe commun dans~$\partial X'$ et qu'il existe un entier~$1\leq i_0\leq n$ tel que~$g_{i_0}$ soit hyperbolique.
Alors l'intersection des~$X'_i$ est un sous-ensemble convexe fermé de la droite géodésique~$X'_{i_0}$.
Il ne contient pas de demi-droite géodésique par hypothèse, donc c'est un segment géodésique~$I$, défini comme l'intersection de deux sous-demi-droites géodésiques~$D_1$ et~$D_2$ de~$X'_{i_0}$.
Comme~$X'$ est un arbre, il existe deux entiers~$i_1$ et~$i_2$ tels que $X'_{i_1}\cap X'_{i_0}\subset D_1$ et $X'_{i_2}\cap X'_{i_0}\subset D_2$.
Pour tout $x'\in\nolinebreak\mathcal{C}$, nous avons vu que $d(x',X'_i)\leq R'/2$ pour tout~$i$.
Si l'on avait $d(x',D_1)>\nolinebreak R'/2$, alors le projeté de~$x'$ sur~$X'_{i_1}$ serait l'extrémité dans~$X'$ de la demi-droite~$D_1$ et l'on aurait $d(x',X'_{i_1})=d(x',D_1)>R'/2$, d'où une contradiction.
Par conséquent, on a $d(x',D_1)\leq\nolinebreak R'/2$ pour tout~$x'\in\mathcal{C}$, et de même $d(x',D_2)\leq R'/2$ pour tout~$x'\in\mathcal{C}$.
Comme~$X'$ est un arbre, on en déduit $d(x',I)\leq R'/2$ pour tout~$x'\in\mathcal{C}$.
Ainsi, le fermé~$\mathcal{C}$ est borné, donc compact.

Supposons que les~$g_i$ admettent exactement un point fixe commun~$\xi'$ dans~$\partial X'$.
Si tous les~$g_i$ étaient elliptiques, il existerait une demi-droite géodésique d'extrémité~$\xi'$ fixée (point par point) par tous les~$g_i$, ce qui contredirait le fait que les~$g_i$ n'ont pas de point fixe commun dans~$X'$.
Par conséquent, il existe un entier~$1\leq i_0\leq n$ tel que~$g_{i_0}$ soit hyperbolique.
L'intersection des~$X'_i$ est un sous-ensemble convexe de la droite géodésique~$X'_{i_0}$ ; comme les~$g_i$ admettent~$\xi'$ comme unique point fixe commun dans~$\partial X'$, cette intersection est une sous-demi-droite géodésique~$D$ de~$X'_{i_0}$ d'extrémité~$\xi'$.
Considérons un entier~$1\leq i\leq n$ tel que~$X'_i\cap X'_{i_0}=D$.
Pour tout~$x'\in\mathcal{C}$, nous avons vu que $d(x',X'_{i_0})\leq R'/2$.
Si l'on avait $d(x',D)>R'/2$, alors le projeté de~$x'$ sur~$X'_i$ serait l'extrémité dans~$X'$ de la demi-droite~$D$ et l'on aurait $d(x',X'_i)=d(x',D)>R'/2$, d'où une contradiction.
Ainsi, on a $d(x',D)\leq R'/2$ pour tout~$x'\in\mathcal{C}$.

Supposons que les~$g_i$ admettent deux points fixes communs~$\xi'_1$ et~$\xi'_2$ \linebreak dans~$\partial X'$.
Par le même argument que ci-dessus, il existe un entier~$1\leq i_0\leq n$ tel que~$g_{i_0}$ soit hyperbolique.
L'axe de translation~$X'_{i_0}$ de~$g_{i_0}$ est la droite~$D$ d'extrémités $\xi'_1$ et~$\xi'_2$.
D'après ce qui précède, on a $d(x',D)\leq R'/2$ pour tout~$x'\in\mathcal{C}$.
\end{dem}

\medskip

Nous sommes maintenant en mesure de démontrer le point~(2) de la proposition~\ref{max sur un ensemble fini}.

\medskip

\noindent
\textbf{Démonstration du point}~(2) \textbf{de la proposition~\ref{max sur un ensemble fini}.}
Si le groupe~$\rho(\Gamma_0)$ admet un point fixe~$x'$ dans~$X'$, alors l'application constante $f : X\rightarrow X'$ d'image~$\{ x'\} $ est $\rho$-équivariante, affine par morceaux, de constante de Lipschitz nulle donc minimale.
\textit{Dans tout le reste de la démonstration, nous supposerons donc $\rho(\Gamma_0)$ sans point fixe dans~$X'$.}
Soit~$\mathcal{D}$ un domaine fondamental connexe pour l'action de~$\Gamma_0$ sur~$X$.
L'ensemble
$$\mathcal{F} = \big\{ \gamma\in\Gamma_0\smallsetminus\{ 1\} ,\quad \gamma\cdot\D\cap\D\neq\emptyset\big\} $$
est un système générateur fini de~$\Gamma_0$.
Comme~$X$ est un arbre, pour tout~$\gamma\in\mathcal{F}$ il existe un unique couple~$(x_{\gamma},y_{\gamma})$ de points de~$\D$ tels que~$\gamma\cdot x_{\gamma}=y_{\gamma}$.
La restriction de~$X$ à~$\D$ induit une bijection entre l'ensemble~$\mathcal{S}$ des applications $\rho$-équivariantes et affines par morceaux $f : X\rightarrow X'$ et l'ensemble~$\mathcal{S}_{\D}$ des applications affines par morceaux $f' : \D\rightarrow X'$ telles que $f'(y_{\gamma})=\rho(\gamma)\cdot f'(x_{\gamma})$ pour tout~$\gamma\in\mathcal{F}$.
Cette bijection préserve les constantes de Lipschitz.
Fixons un réel~$C'>C$ et notons~$\mathcal{S}_{\D,C'}$ l'ensemble des éléments de~$\mathcal{S}_{\D}$ de constante de Lipschitz~$\leq C'$ ; il est fermé et uniformément équicontinu dans l'ensemble des applications continues de~$\D$ dans~$X'$, muni de la topologie de la convergence uniforme.
Pour pouvoir appliquer le théorème d'Ascoli, nous allons nous ramener à un sous-ensemble de~$\mathcal{S}_{\D,C'}$ à valeurs dans un compact de~$X'$, en distinguant plusieurs cas selon le nombre de points fixes de~$\rho(\Gamma_0)$ dans~$\partial X'$.
Commençons par remarquer que si l'on pose
$$R\ =\ \max\big\{ d(x,\gamma\cdot x),\quad x\in\D,\ \gamma\in\mathcal{F}\big\} >0,$$
alors pour tout~$f'\in\mathcal{S}_{\D,C'}$ et tout~$x\in\D$ on a $d(f'(x),\rho(\gamma)\cdot f'(x))\leq C'R$ pour tout~$\gamma\in\mathcal{F}$.

Si $\rho(\Gamma_0)$ n'a pas de point fixe dans~$\partial X'$, alors~$\D'$ est compact d'après le lemme~\ref{lemme axes de translation et points fixes}.
On peut donc appliquer directement le théorème d'Ascoli, qui affirme que l'ensemble~$\mathcal{S}_{\D,C'}$ est compact pour la topologie de la convergence uniforme.
On en déduit l'existence d'un élément $f'\in\mathcal{S}_{\D,C'}$, et donc d'un élément~$f\in\mathcal{S}$, de constante de Lipschitz~$C$ minimale.

Supposons que $\rho(\Gamma_0)$ admette un unique point fixe~$\xi'\in\partial X'$.
D'après le lemme~\ref{lemme axes de translation et points fixes}, il existe une demi-droite géodésique~$D$ de~$X'$ d'extrémité~$\xi'$ telle que $d(f'(x),D)\leq C'R/2$ pour tout~$f'\in\mathcal{S}_{\D,C'}$ et tout~$x\in\D$.
Fixons une isométrie hyperbolique~$g\in\Isom(X')$ dont l'axe de translation contient~$D$ et dont~$\xi'$ est le point fixe répulsif dans~$\partial X'$.
Soit~$(f'_n)_{n\in\N}$ une suite d'éléments de~$\mathcal{S}_{\D,C'}$ dont la constante de Lipschitz tend vers~$C$.
Si l'on note~$x'_0$ l'extrémité de~$D$ dans~$X'$, alors pour tout~$n\in\N$ il existe un entier~$m_n\in\N$ tel que~$g^{m_n}\cdot f'_n(\D)$ soit inclus dans le compact
$$\D' = \Big\{ x'\in X',\quad d(x',x'_0) \leq \lambda(g)+\frac{C'R}{2}\Big\} .$$
Pour tout~$\gamma\in\Gamma_0$, la suite $(g^{m_n}\rho(\gamma)g^{-m_n})_{n\in\N}$ est bornée : en effet, comme $\rho(\gamma)$ fixe~$\xi'\in\partial X'$ par hypothèse, il existe une demi-droite~$D_{\gamma}$ incluse dans~$D$ telle que~$\rho(\gamma)\cdot D_{\gamma}\subset D$ ; pour tout~$x'\in D_{\gamma}$ on a $g^{m_n}\rho(\gamma)g^{-m_n}\cdot x'=\rho(\gamma)\cdot x'$ pour tout~$n\in\N$, donc~$(g^{m_n}\rho(\gamma)g^{-m_n})_{n\in\N}$ est bornée par définition de la topologie sur~$\Isom(X')$ (\textit{cf.} paragraphe~\ref{Action d'un groupe sur un arbre}).
Quitte à extraire une sous-suite, on peut donc supposer que pour tout~$\gamma\in\mathcal{F}$ la suite~$(g^{m_n}\rho(\gamma)g^{-m_n})_{n\in\N}$ converge.
En particulier, il existe un entier~$N\in\N$ tel que pour tout~$n\geq N$ et tout~$\gamma\in\mathcal{F}$, les isométries $g^{m_n}\rho(\gamma)g^{-m_n}$ et~$g^{m_N}\rho(\gamma)g^{-m_N}$ coïncident sur le compact~$\D'$.
Par conséquent, pour tout $n\geq N$ et tout~$\gamma\in\mathcal{F}$ on a
\begin{eqnarray*}
\big(g^{m_n-m_N}\cdot f'_n\big)(y_{\gamma}) & = & g^{m_n-m_N}\rho(\gamma)\cdot f'_n(x_{\gamma})\\
& = & g^{-m_N} \big(g^{m_n}\rho(\gamma) g^{-m_n}\big) \cdot \big(g^{m_n} \cdot f'_n(x_{\gamma})\big)\\
& = & g^{-m_N} \big(g^{m_N}\rho(\gamma) g^{-m_N}\big) \cdot \big(g^{m_n} \cdot f'_n(x_{\gamma})\big)\\
& = & \rho(\gamma) \cdot \big(g^{m_n-m_N}\cdot f'_n\big)(x_{\gamma}).
\end{eqnarray*}
Comme de plus $g$ est une isométrie de~$X'$, l'application $g^{m_n-m_N}\cdot f'_n$ a même constante de Lipschitz que~$f'_n$ ; elle appartient donc à~$\mathcal{S}_{\D,C'}$ pour tout~$n\geq N$.
L'image de~$g^{m_n-m_N}\cdot f'_n$ étant incluse dans le compact~$g^{-m_N}\cdot\D'$, on peut appliquer le théorème d'Ascoli, qui affirme que la suite~$(g^{m_n-m_N}\cdot f'_n)_{n\in\N}$ est d'adhérence compacte dans le fermé~$\mathcal{S}_{\D,C'}$ pour la topologie de la convergence uniforme.
En particulier, cette suite admet une sous-suite qui converge vers un certain~$f'\in\mathcal{S}_{\D,C'}$.
Comme la constante de Lipschitz de~$f'_n$ tend vers~$C$, celle de~$g^{m_n-m_N}\cdot f'_n$ aussi, et $f'$ est de constante de Lipschitz~$C$ minimale.

Le cas où~$\rho(\Gamma_0)$ admet deux points fixes~$\xi'_1$ et~$\xi'_2$ dans~$\partial X'$ se traite de manière analogue.
Plus précisément, soit~$D$ la droite géodésique de~$X'$ d'extrémités~$\xi'_1$ et~$\xi'_2$.
D'après le lemme~\ref{lemme axes de translation et points fixes}, on a $d(f'(x),D)\leq C'R/2$ pour tout~$f'\in\mathcal{S}_{\D,C'}$ et tout~$x\in\D$.
Fixons une isométrie hyperbolique $g\in\Isom(X')$ d'axe de translation~$D$.
Soit~$(f'_n)_{n\in\N}$ une suite d'éléments de~$\mathcal{S}_{\D,C'}$ dont la constante de Lipschitz tend vers~$C$.
Si~$x'_0$ désigne un point quelconque fixé de~$D$, alors pour tout~$n\in\N$ il existe un entier~$m_n\in\Z$ tel que~$g^{m_n}\cdot f'_n(\D)$ soit inclus dans le compact
$$\D' = \Big\{ x'\in X',\quad d(x',x'_0)\leq\lambda(g)+\frac{C'R}{2}\Big\} .$$
Pour tout~$\gamma\in\Gamma_0$, la suite~$(g^{m_n}\rho(\gamma)g^{-m_n})_{n\in\N}$ est bornée : en effet, la droite~$D$ est (globalement) stable par~$\rho(\gamma)$ car $\rho(\gamma)$ fixe ses extrémités $\xi'_1,\xi'_2\in\partial X'$ par hypothèse ; pour tout~$x'\in D$ on a $g^{m_n}\rho(\gamma)g^{-m_n}\cdot x'=\rho(\gamma)\cdot x'$ pour tout~$n\in\N$, donc~$(g^{m_n}\rho(\gamma)g^{-m_n})_{n\in\N}$ est bornée par définition de la topologie sur~$\Isom(X')$ (\textit{cf.} paragraphe~\ref{Action d'un groupe sur un arbre}).
On conclut comme dans le cas où~$\rho(\Gamma_0)$ n'admet qu'un seul point fixe dans~$\partial X'$.
\hfill\qedsymbol

\bigskip

Notons que~$C\geq 0$ est la constante de Lipschitz minimale de toutes les applications lipschitziennes et $\rho$-équivariantes de~$X$ dans~$X'$ (non nécessairement affines par morceaux).
En effet, si $f : X\rightarrow X'$ est $\rho$-équivariante et $C$-lipschitzienne, alors l'application affine par morceaux $f_a : X\rightarrow X'$ telle que $f_a(s)=f(s)$ pour tout sommet~$s$ de~$X$ (donnée par la remarque~\ref{remarque affine par morceaux}) est encore $\rho$-équivariante et $C$-lipschitzienne.

%%%%%%%%%%%%%%%%%%%%%%%%%
\subsection{Arêtes $f$-maximales}

Fixons désormais deux arbres réels simpliciaux $X$ et~$X'$ de valence~$\geq 2$, un sous-groupe discret sans torsion~$\Gamma_0$ de~$\Isom(X)$ tel que le graphe~$\Gamma_0\backslash X$ soit fini et un morphisme de groupes $\rho : \Gamma_0\rightarrow G$.
Soit $C\geq 0$ la constante de Lipschitz minimale donnée par le point~(2) de la proposition~\ref{max sur un ensemble fini}.
Fixons une application $\rho$-équivariante, affine par morceaux et $C$-lipschitzienne $f : X\rightarrow X'$.

Nous dirons qu'une arête~$e$ de $X$ est $f$\textit{-maximale} si la restriction de~$f$ à~$e$ est affine de constante~$C$.
Comme $f$ est $\rho$-équivariante et comme $\Gamma_0$ et~$\rho(\Gamma_0)$ sont des groupes d'isométries, une arête~$e$ de~$X$ est $f$-maximale si et seulement l'arête~$\gamma\cdot\nolinebreak e$ est $f$-maximale pour tout~$\gamma\in\Gamma_0$.
En particulier, la notion d'arête maximale passe au quotient~$\Gamma_0\backslash X$.

Pour toutes arêtes~$e_1$ et~$e_2$ de~$X$ incidentes en un sommet~$s$, notons $e_1\sim_f e_2$ si $f(e_1)\cap f(e_2)$ est d'intérieur non vide dans~$X'$.
Ceci définit une relation d'équivalence~$\sim_f$ sur l'ensemble des arêtes de~$X$ d'extrémité~$s$.
Comme précédemment, cette relation induit une relation d'équivalence~$\sim_f$ sur l'ensemble des arêtes de~$\Gamma_0\backslash X$ incidentes en~$\overline{s}$, où $\overline{s}$ désigne l'image de~$s$ dans~$\Gamma_0\backslash X$.

Pour démontrer le point~(3) de la proposition~\ref{max sur un ensemble fini}, nous aimerions prouver l'existence d'un lacet de~$\Gamma_0\backslash X$ passant au plus deux fois par chaque sommet de~$\Gamma_0\backslash X$ et tel que la restriction de~$f$ aux relevés dans~$X$ de ce lacet soit affine de constante~$C$.
Cette dernière condition se traduit par le fait que le lacet est entièrement formé d'arêtes $f$-maximales et que deux arêtes successives appartiennent toujours à des classes d'équivalence différentes pour la relation~$\sim_f$.
L'existence d'un tel lacet est assurée si~$f$ vérifie une certaine condition de minimalité : il suffit que l'ensemble des arêtes $f$-maximales de~$X$ soit minimal.
Avant de démontrer cette existence (lemme~\ref{boucle d'aretes f-maximales}), commençons par établir qu'une condition nécessaire est vérifiée.

\begin{lem}\label{aretes f-maximales}
Supposons l'ensemble des arêtes $f$-maximales de~$X$ minimal.
En tout sommet de~$X$ qui est extrémité d'une arête $f$-maximale, il existe au moins deux classes d'équivalence d'arêtes $f$-maximales pour la relation~$\sim_f$.
\end{lem}

\begin{dem}
Soit~$s_0$ un sommet de~$X$.
Supposons par l'absurde que toutes les arêtes $f$-maximales d'extrémité~$s_0$ appartiennent à la même classe d'équivalence pour la relation~$\sim_f$ : cela signifie qu'il existe un segment géodésique~$I$ d'intérieur non vide dans~$X'$, d'extrémité~$f(s_0)$, qui est contenu dans l'image~$f(e)$ de toute arête $f$-maximale~$e$ d'extrémité~$s_0$.
D'après la remarque~\ref{remarque affine par morceaux}, pour tout~$x'\in I$ il existe une unique application affine par morceaux~$f_{x'} : X\rightarrow X'$ telle que $f_{x'}(\gamma\cdot s_0)=\rho(\gamma)\cdot x'$ pour tout~$\gamma\in\Gamma_0$ et $f_{x'}(s)=f(s)$ pour tout sommet~$s\notin\Gamma_0\cdot s_0$ ; elle est $\rho$-équivariante.
Montrons que pour~$x'\in I\smallsetminus\{ f(s_0)\} $ suffisamment proche de~$f(s_0)$, l'application~$f_{x'}$ est $C$-lipschitzienne et l'ensemble des arêtes $f_{x'}$-maximales est strictement inclus dans l'ensemble des arêtes $f$-maximales, ce qui contredira la minimalité de~$C$ et de l'ensemble des arêtes $f$-maximales.

Si~$e$ est une arête de~$X$ dont les extrémités n'appartiennent pas à l'orbite~$\Gamma_0\cdot s_0$, alors la restriction de~$f_{x'}$ à~$e$ est égale à la restriction de~$f$ à~$e$.
En particulier, $e$ est $f_{x'}$-maximale si et seulement si elle est $f$-maximale.

Soit~$e$ une arête de~$X$ d'extrémités~$s_0$ et~$s\notin\Gamma_0\cdot s_0$.
La restriction de~$f$ (resp. de~$f_{x'}$) à~$e$ est affine de constante
$$\frac{d(f(s_0),f(s))}{d(s_0,s)} \quad\quad\mathrm{\Big(resp.}\quad \frac{d(x',f(s))}{d(s_0,s)}\mathrm{\Big)}.$$
Si~$e$ n'est pas $f$-maximale, alors $e$ n'est pas $f_{x'}$-maximale pour $x'$ suffisamment proche de~$f(s_0)$, par continuité de l'application $x'\mapsto d(x',f(s))$.
Si~$e$ est $f$-maximale, alors pour~$x'\in I\smallsetminus\{ f(s_0)\} $ on a $d(x',f(s))<d(f(s_0),f(s))$ donc $e$ n'est pas $f_{x'}$-maximale.

Soit~$e$ une arête de~$X$ d'extrémités~$s_0$ et~$\gamma\cdot s_0$, où~$\gamma\in\Gamma_0\smallsetminus\{ 1\} $.
La restriction de~$f$ (resp. de~$f_{x'}$) à~$e$ est affine de constante
$$\frac{d(f(s_0),\rho(\gamma)\cdot f(s_0))}{d(s_0,\gamma\cdot s_0)} \quad\quad\mathrm{\Big(resp.}\quad \frac{d(x',\rho(\gamma)\cdot x')}{d(s_0,\gamma\cdot s_0)}\mathrm{\Big)}.$$
Si~$e$ n'est pas $f$-maximale, alors $e$ n'est pas $f_{x'}$-maximale pour $x'$ suffisamment proche de~$f(s_0)$, par continuité de l'application $x'\mapsto d(x',\rho(\gamma)\cdot x')$.
Supposons que $e$ est $f$-maximale et montrons que $e$ n'est pas $f_{x'}$-maximale pour $x'\in I\smallsetminus\{ f(s_0)\} $.
Comme~$e$ est $f$-maximale, $\gamma^{-1}\cdot e$ l'est aussi, donc l'intersection
$$f(\gamma^{-1}\cdot e)\cap f(e)\ =\ [\rho(\gamma)^{-1}\cdot f(s_0),f(s_0)]\cap [f(s_0),\rho(\gamma)\cdot f(s_0)]$$
contient~$I$.
Si~$\rho(\gamma)$ est hyperbolique, alors pour $x'\in I\smallsetminus\{ f(s_0)\} $ on a, d'après~(\ref{mu en fonction de lambda}),
\begin{eqnarray*}
d(x',\rho(\gamma)\cdot x') & = & \lambda(\rho(\gamma)) + 2d(x',\mathcal{A}_{\rho(\gamma)})\\
& < & \lambda(\rho(\gamma)) + 2d(f(s_0),\mathcal{A}_{\rho(\gamma)})\\
& = & d(f(s_0),\rho(\gamma)\cdot f(s_0)),
\end{eqnarray*}
donc $e$ n'est pas $f_{x'}$-maximale.
Si~$\rho(\gamma)$ est elliptique et si~$X'_{\rho(\gamma)}$ désigne l'ensemble de ses points fixes dans~$X'$, alors pour $x'\in I\smallsetminus\{ f(s_0)\} $ on~a, d'après~(\ref{distance de x a gx pour g non hyperbolique}),
$$d(x',\rho(\gamma)\cdot x') = 2d(x',X'_{\rho(\gamma)}) < 2d(f(s_0),X'_{\rho(\gamma)}) = d(f(s_0),\rho(\gamma)\cdot f(s_0)),$$
donc $e$ n'est pas $f_{x'}$-maximale, ce qui termine la démonstration du lemme~\ref{aretes f-maximales}.
\end{dem}

%%%%%%%%%%%%%%%%%%%%%%%%%
\subsection{Un lacet d'étirement maximal}

Démontrons à présent l'existence d'un lacet de~$\Gamma_0\backslash X$ passant au plus deux fois par chaque sommet de~$\Gamma_0\backslash X$ et tel que la restriction de~$f$ aux relevés dans~$X$ de ce lacet soit affine de constante~$C$.

\begin{lem}\label{boucle d'aretes f-maximales}
Supposons l'ensemble des arêtes $f$-maximales de~$X$ minimal.
Il existe un élément~$\gamma\in\Gamma_0\smallsetminus\{ 1\} $ dont l'image dans~$\Gamma_0\backslash X$ de l'axe de translation~$\mathcal{A}_{\gamma}$ est un lacet passant au plus deux fois par chaque sommet de~$\Gamma_0\backslash X$, dont toutes les arêtes sont $f$-maximales et dont deux arêtes consécutives appartiennent toujours à des classes d'équivalence différentes pour la relation~$\sim_f$.
\end{lem}

\begin{dem}
Soit $e_0$ une arête $f$-maximale de~$X$, d'extrémités~$s_0$ et~$s_1$.
Le lemme~\ref{aretes f-maximales} permet de construire par récurrence une suite~$(s_n)_{n\in\N}$ de sommets de~$X$ telle que pour tout~$n\in\N$, les sommets~$s_n$ et~$s_{n+1}$ soient adjacents, l'arête~$e_n$ d'extrémités~$s_n$ et~$s_{n+1}$ soit $f$-maximale, et les arêtes~$e_n$ et~$e_{n+1}$ appartiennent à des classes d'équivalence différentes pour la relation~$\sim_f$.
Comme~$X$ est un arbre, les sommets~$s_n$ sont deux à deux distincts et pour tous~$m<n$, la réunion des arêtes~$e_i$, où~$m\leq i\leq n$, est le segment géodésique d'extrémités~$s_m$ et~$s_n$, de longueur~$n-m+1$.
Soit~$N\geq 1$ le nombre de sommets du graphe fini~$\Gamma_0\backslash X$.
D'après le principe des tiroirs de Dirichlet, il existe trois entiers
$$0\leq n_1<n_2<n_3\leq 2N$$
tels que les sommets $s_{n_1}$, $s_{n_2}$ et~$s_{n_3}$ appartiennent à la même orbite de~$\Gamma_0$.
Choisissons-les de sorte que l'entier $n_3-n_1$ soit minimal, ce qui assure que l'image du segment géodésique $[s_{n_1},s_{n_3}]$ dans~$\Gamma_0\backslash X$ est un lacet passant au plus deux fois par chaque sommet de~$\Gamma_0\backslash X$.
Soient $\gamma_1,\gamma_2\in\Gamma_0\smallsetminus\{ 1\} $ tels que~$s_{n_2}=\gamma_1\cdot s_{n_1}$ et $s_{n_3}=\gamma_2\cdot s_{n_2}$.
La figure~4 illustre notre raisonnement.
Pour toute arête~$e$ de~$X$, nous notons~$\overline{e}$ son image dans~$\Gamma_0\backslash X$.

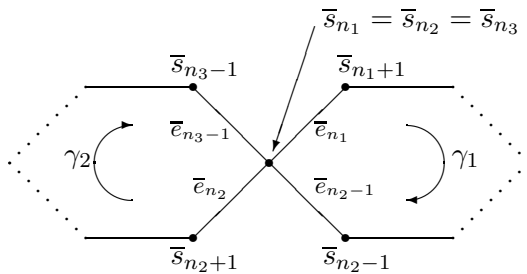
\begin{figure}
\begin{center}\setlength{\unitlength}{1mm}
\begin{picture}(70,34)(-35,-14)
\put(0,0){\line(1,1){10}}
\put(10,10){\line(1,0){14.14}}
\multiput(24.14,10)(1.25,-1.25){9}{\circle*{0.4}}
\put(0,0){\line(1,-1){10}}
\put(10,-10){\line(1,0){14.14}}
\multiput(24.14,-10)(1.25,1.25){9}{\circle*{0.4}}

\put(0,0){\line(-1,1){10}}
\put(-10,10){\line(-1,0){14.14}}
\multiput(-24.14,10)(-1.25,-1.25){9}{\circle*{0.4}}
\put(0,0){\line(-1,-1){10}}
\put(-10,-10){\line(-1,0){14.14}}
\multiput(-24.14,-10)(-1.25,1.25){9}{\circle*{0.4}}

\put(0,0){\circle*{1}}
\put(10,10){\circle*{1}}
\put(10,-10){\circle*{1}}
\put(-10,10){\circle*{1}}
\put(-10,-10){\circle*{1}}

\put(6,3.5){\footnotesize $\overline{e}_{n_1}$}
\put(6,-4){\footnotesize $\overline{e}_{n_2-1}$}
\put(-13,3.5){\footnotesize $\overline{e}_{n_3-1}$}
\put(-10,-4){\footnotesize $\overline{e}_{n_2}$}

\put(9,12){$\overline{s}_{n_1+1}$}
\put(7,-13.5){$\overline{s}_{n_2-1}$}
\put(-13,12){$\overline{s}_{n_3-1}$}
\put(-13,-13.5){$\overline{s}_{n_2+1}$}

\put(18,0){\oval(10,10)[r]}
\put(19,-5){\vector(-1,0){1}}
\put(24,0){$\gamma_1$}
\put(-18,0){\oval(10,10)[l]}
\put(-19,5){\vector(1,0){1}}
\put(-27,0){$\gamma_2$}

\put(6,18){\vector(-1,-3){5.5}}
\put(7,18){$\overline{s}_{n_1}=\overline{s}_{n_2}=\overline{s}_{n_3}$}
\end{picture}
\end{center}
\caption{Illustration du lemme~\ref{boucle d'aretes f-maximales}, dans~$\Gamma_0\backslash X$}
\end{figure}

Supposons que les arêtes~$\overline{e}_{n_2-1}$ et~$\overline{e}_{n_1}$ appartiennent à des classes d'équivalence différentes pour la relation~$\sim_f$.
Le projeté de~$s_{n_2}$ sur~$\mathcal{A}_{\gamma_1}$ appartient à la fois aux segments géodésiques~$[s_{n_1},s_{n_2}]$ et $[s_{n_2},\gamma_1\cdot s_{n_2}]$ ; comme $e_{n_2-1}\not\sim_f\gamma_1\cdot e_{n_1}$ par hypothèse, on a~$s_{n_2}\in\mathcal{A}_{\gamma_1}$.
Ainsi, l'image de~$\mathcal{A}_{\gamma_1}$ dans~$\Gamma_0\backslash X$ est l'union des~$\overline{e}_i$ où $n_1\leq i\leq n_2-1$ ; c'est un lacet simple de~$\Gamma_0\backslash X$ dont toutes les arêtes sont $f$-maximales et dont deux arêtes consécutives appartiennent toujours à des classes d'équivalence différentes pour la relation~$\sim_f$.

De même, si les arêtes~$\overline{e}_{n_3-1}$ et~$\overline{e}_{n_2}$ appartiennent à des classes d'équivalence différentes pour la relation~$\sim_f$, alors l'image de~$\mathcal{A}_{\gamma_2}$ dans~$\Gamma_0\backslash X$ est l'union des~$\overline{e}_i$ où $n_2\leq i\leq n_3-1$ ; c'est un lacet simple de~$\Gamma_0\backslash X$ dont toutes les arêtes sont $f$-maximales et dont deux arêtes consécutives appartiennent toujours à des classes d'équivalence différentes pour la relation~$\sim_f$.

Supposons~$\overline{e}_{n_2-1}\sim_f\overline{e}_{n_1}$ et~$\overline{e}_{n_3-1}\sim_f\overline{e}_{n_2}$.
Comme par construction les arêtes $\overline{e}_{n_2-1}$ et~$\overline{e}_{n_2}$ appartiennent à des classes d'équivalence différentes pour~$\sim_f$, il en est de même des arêtes $\overline{e}_{n_3-1}$ et~$\overline{e}_{n_1}$.
Comme précédemment, l'image de~$\mathcal{A}_{\gamma_2\gamma_1}$ dans~$\Gamma_0\backslash X$ est l'union des~$\overline{e}_i$ où $n_1\leq i\leq n_3-1$ ; c'est la concaténation de deux lacets simples de~$\Gamma_0\backslash X$, dont toutes les arêtes sont $f$-maximales et dont deux arêtes consécutives appartiennent toujours à des classes d'équivalence différentes pour la relation~$\sim_f$.
Ceci achève la démonstration du lemme~\ref{boucle d'aretes f-maximales}.
\end{dem}

\medskip

Nous pouvons à présent démontrer le point~(3) de la proposition~\ref{max sur un ensemble fini}.

\medskip

\noindent
\textbf{Démonstration du point}~(3) \textbf{de la proposition~\ref{max sur un ensemble fini}.}
Soit $f : X\rightarrow\nolinebreak X'$ une application $\rho$-équivariante, affine par morceaux et $C$-lipschitzienne telle que l'ensemble des arêtes $f$-maximales de~$X$ soit minimal.
D'après le lemme~\ref{boucle d'aretes f-maximales}, il existe un élément~$\gamma\in\Gamma_0\smallsetminus\{ 1\} $ dont l'image dans~$\Gamma_0\backslash X$ de l'axe de translation~$\mathcal{A}_{\gamma}$ est un lacet passant au plus deux fois par chaque sommet de~$\Gamma_0\backslash X$, dont toutes les arêtes sont $f$-maximales et dont deux arêtes consécutives appartiennent toujours à des classes d'équivalence différentes pour la relation~$\sim_f$.
Quitte à remplacer~$\gamma$ par un conjugué dans~$\Gamma_0$, on peut supposer que~$\mathcal{A}_{\gamma}$ rencontre~$\mathcal{D}$ ; d'après~(\ref{mu en fonction de lambda}) on a alors
$$d(x_0,\gamma\cdot x_0)\, =\, \lambda(\gamma) + 2d(x_0,\mathcal{A}_{\gamma})\, \leq\, 4N,$$
\textit{i.e.}~$\gamma\in F$.
Montrons que
$$\frac{\lambda\big(\rho(\gamma)\big)}{\lambda(\gamma)}\, =\, C,$$
ce qui suffira, d'après~(\ref{inegalite sup-inf}), à établir le point~(3) de la proposition~\ref{max sur un ensemble fini}.
Nous distinguons les cas où~$\rho(\gamma)$ est hyperbolique ou elliptique.

Supposons~$\rho(\gamma)$ hyperbolique et soit~$s$ un sommet de~$\mathcal{A}_{\gamma}$.
Le projeté de~$f(s)$ sur~$\mathcal{A}_{\rho(\gamma)}$ appartient à la fois aux segments géodésiques $[f(s),\rho(\gamma)^{-1}\cdot f(s)]=[f(s),f(\gamma^{-1}\cdot s)]$ et $[f(s),\rho(\gamma)\cdot f(s)]=[f(s),f(\gamma\cdot s)]$.
Comme deux arêtes consécutives de~$\mathcal{A}_{\gamma}$ appartiennent toujours à des classes d'équivalence différentes pour la relation~$\sim_f$, on en déduit $f(s)\in\mathcal{A}_{\rho(\gamma)}$.
Ainsi, on a
$$\frac{d(f(s),f(\gamma\cdot s))}{d(s,\gamma\cdot s)}\ =\ \frac{d(f(s),\rho(\gamma)\cdot f(s))}{d(s,\gamma\cdot s)}\ =\ \frac{\lambda(\rho(\gamma))}{\lambda(\gamma)}.$$
Comme toutes les arêtes de~$\mathcal{A}_{\gamma}$ sont $f$-maximales, ce quotient est égal à~$C$.

Supposons~$\rho(\gamma)$ elliptique, de sorte que~$\lambda(\rho(\gamma))=0$, et soit~$s$ un sommet de~$\mathcal{A}_{\gamma}$.
Le projeté de~$f(s)$ sur l'ensemble des points fixes de~$\rho(\gamma)$ appartient à la fois aux segments géodésiques $[f(s),\rho(\gamma)^{-1}\cdot f(s)]=[f(s),f(\gamma^{-1}\cdot s)]$ et $[f(s),\rho(\gamma)\cdot f(s)]=[f(s),f(\gamma\cdot s)]$.
Comme deux arêtes consécutives de~$\mathcal{A}_{\gamma}$ appartiennent toujours à des classes d'équivalence différentes pour la relation~$\sim_f$, on en déduit que $f(s)$ est un point fixe de~$\rho(\gamma)$.
Ainsi, on a
$$f(\gamma\cdot s)\, =\, \rho(\gamma)\cdot f(s)\, =\, f(s),$$
donc $C=0$, ce qui termine la démonstration.
\hfill\qedsymbol

%%%%%%%%%%%%%%%%%%%%%%%%%%%%%%%%%%%%%%%%%%%%%%%%%%%%%%%%%%%%
\section{Déformation des quotients compacts de $(G\times G)/\Delta_G$}\label{Deformation des quotients compacts}

Dans cette partie, nous déduisons les théorèmes~\ref{condition d'admissibilite} puis~\ref{ouvert} de la proposition~\ref{max sur un ensemble fini} et d'une propriété des déformations de réseaux cocompacts sans torsion de~$G$ (lemme~\ref{Reste un reseau cocompact sans torsion}).

%%%%%%%%%%%%%%%%%%%%%%%%%
\subsection{Démonstration du théorème~\ref{condition d'admissibilite}}\label{Demonstration du theoreme sur l'admissibilite}

Soient $X$ et~$X'$ deux arbres réels simpliciaux de valence~$\geq 2$ et $\Gamma_0$ un sous-groupe discret sans torsion de $\Isom(X)$ tel que le graphe~$\Gamma_0\backslash X$ soit fini.
Fixons des points~$x_0\in X$ et~$x'_0\in\nolinebreak X'$.
Nous dirons qu'un morphisme de groupes $\rho : \Gamma_0\rightarrow\Isom(X')$ est \textit{admissible} si pour tout~$R>0$ on a $d(x'_0,\rho(\gamma)\cdot x'_0)\leq d(x_0,\gamma\cdot x_0)-R$ pour presque tout~$\gamma\in\Gamma_0$.
Cette condition ne dépend pas du choix de~$x_0$ et~$x'_0$.
La proposition~\ref{max sur un ensemble fini} implique le résultat suivant, dont le théorème~\ref{condition d'admissibilite} est un cas particulier.

\begin{theo}\label{condition d'admissibilite sur les arbres}
Soient $X$ et~$X'$ deux arbres réels simpliciaux de valence~$\geq 2$ et $\Gamma_0$ un sous-groupe discret sans torsion de~$\Isom(X)$ tel que le graphe~$\Gamma_0\backslash X$ soit fini.
Fixons des points~$x_0\in X$ et~$x'_0\in X'$.
Pour tout morphisme de groupes $\rho : \Gamma_0\rightarrow\Isom(X')$, notons~$C_{\rho}\geq 0$ la borne inférieure des réels~$t\geq 0$ pour lesquels il existe~$t'\geq 0$ tel que $d(x'_0,\rho(\gamma)\cdot x'_0)\leq td(x_0,\gamma\cdot x_0)+t'$ pour tout~$\gamma\in\Gamma_0$.
\begin{enumerate}
	\item Soit $F$ le sous-ensemble fini de~$\Gamma_0\smallsetminus\{ 1\} $ formé des éléments~$\gamma$ tels que $d(x_0,\gamma\cdot x_0)\leq\nolinebreak 4N$, où $N$ désigne le nombre de sommets de~$\Gamma_0\backslash X$.
	Pour tout morphisme $\rho : \Gamma_0\rightarrow\Isom(X')$ on~a
	$$C_{\rho}\ =\ \sup_{\gamma\in\Gamma_0\smallsetminus\{ 1\} } \frac{\lambda(\rho(\gamma))}{\lambda(\gamma)}\ =\ \max_{\gamma\in F}\, \frac{\lambda(\rho(\gamma))}{\lambda(\gamma)}.$$
	\item Un morphisme $\rho : \Gamma_0\rightarrow\Isom(X')$ est admissible si et seulement si~$C_{\rho}<1$.
\end{enumerate}
\end{theo}

\medskip

Comme dans la proposition~\ref{max sur un ensemble fini}, l'hypothèse que $\Gamma_0$ est un sous-groupe discret sans torsion de~$\Isom(X)$ tel que le graphe~$\Gamma_0\backslash X$ soit fini est équivalente à l'hypothèse que $\Gamma_0$ est un groupe libre de type fini agissant sur~$X$ librement, minimalement, par isométries (\textit{cf.} paragraphe~\ref{groupes discrets d'isometries}).

\medskip

\begin{dem}
Soit $\rho : \Gamma_0\rightarrow\Isom(X')$ un morphisme de groupes.
D'après la proposition~\ref{max sur un ensemble fini}, il existe une application $\rho$-équivariante affine par morceaux $f : X\rightarrow X'$ de constante de Lipschitz~$C$ minimale.
Commençons par remarquer que pour tout $\gamma\in\Gamma_0\smallsetminus\{ 1\} $ on a
$$\frac{\lambda(\rho(\gamma))}{\lambda(\gamma)}\ \leq\ C_{\rho}.$$
En effet, c'est évident si~$\rho(\gamma)$ est elliptique, et si~$\rho(\gamma)$ est hyperbolique cela résulte du fait que la suite
\begin{equation}\label{mu circ rho et mu}
\Big(d(x'_0,\rho(\gamma)^n\cdot x'_0) - \frac{\lambda(\rho(\gamma))}{\lambda(\gamma)}\,d(x_0,\gamma^n\cdot x_0)\Big)_{n\geq 1}
\end{equation}
est constante d'après~(\ref{mu en fonction de lambda}).
D'autre part, pour tout~$\gamma\in\Gamma_0$ on a
\begin{eqnarray*}
d(x'_0,\rho(\gamma)\cdot x'_0) & \leq & d\big(x'_0,f(x_0)\big) + d\big(f(x_0),\rho(\gamma)\cdot f(x_0)\big) + d\big(\rho(\gamma)\cdot f(x_0),\rho(\gamma)\cdot x'_0\big)\\
& = & d\big(f(x_0),f(\gamma\cdot x_0)\big) + 2\,d\big(x'_0,f(x_0)\big)\\
& \leq & C\,d(x_0,\gamma\cdot x_0) + 2\,d\big(x'_0,f(x_0)\big),
\end{eqnarray*}
donc~$C_{\rho}\leq C$.
Ainsi, on a
$$\sup_{\gamma\in\Gamma_0\smallsetminus\{ 1\} } \frac{\lambda(\rho(\gamma))}{\lambda(\gamma)}\ \leq\ C_{\rho}\ \leq\ C.$$
D'après la proposition~\ref{max sur un ensemble fini}, il existe un élément~$\gamma'\in F$ tel que
\begin{equation}\label{egalite pour un element de F}
\frac{\lambda(\rho(\gamma'))}{\lambda(\gamma')}\ =\ \sup_{\gamma\in\Gamma_0\smallsetminus\{ 1\} } \frac{\lambda(\rho(\gamma))}{\lambda(\gamma)}\ =\ C_{\rho}\ =\ C,
\end{equation}
ce qui prouve le point~(1) du théorème~\ref{condition d'admissibilite sur les arbres}.
Comme l'application~$\mu$ est propre et le groupe~$\Gamma_0$ discret dans~$G$, si~$C_{\rho}<1$ alors $\rho$ est admissible, et si~$\rho$ est admissible alors~$C_{\rho}\leq 1$.
D'après~(\ref{egalite pour un element de F}), si~$C_{\rho}=1$ alors $\lambda(\rho(\gamma'))=\lambda(\gamma')$ ; comme la suite donnée par~(\ref{mu circ rho et mu}) est constante pour~$\gamma=\gamma'$, le morphisme~$\rho$ n'est pas admissible.
Ceci achève la démonstration du point~(2) du théorème~\ref{condition d'admissibilite sur les arbres}.
\end{dem}

\medskip

Le théorème~\ref{condition d'admissibilite} s'obtient à partir du théorème~\ref{condition d'admissibilite sur les arbres} en prenant pour~$X$ et~$X'$ l'arbre de Bruhat-Tits de~$G$ et pour~$x_0$ et~$x'_0$ le point donné par~(\ref{mu distance}).

%%%%%%%%%%%%%%%%%%%%%%%%%
\subsection{Déformation des réseaux cocompacts sans torsion de $G$}\label{Deformation}

Soient $\kkk$ un corps local ultramétrique et $G$ l'ensemble des $\kkk$-points d'un $\kkk$-groupe algébrique semi-simple connexe de $\kkk$-rang~un.
Il est bien connu que pour tout réseau cocompact sans torsion~$\Gamma_0$ de~$G$, il existe un voisinage~$\mathcal{V}_0$ de l'inclusion canonique dans~$\Hom(\Gamma_0,G)$ tel que pour tout~$\sigma\in\mathcal{V}_0$, le groupe~$\sigma(\Gamma_0)$ soit un réseau cocompact sans torsion de~$G$ (\cite{lub}, prop.~1.7 et th.~2.1).
De plus, comme les applications~$\lambda$ et~$\mu$ sont continues à valeurs discrètes, pour toute partie finie~$F$ de~$\Gamma_0$ il existe un voisinage~$\mathcal{V}_F$ de l'inclusion canonique dans~$\Hom(\Gamma_0,G)$ tel que pour tout~$\sigma\in\mathcal{V}_F$ on ait $\lambda(\sigma(\gamma))=\lambda(\gamma)$ et $\mu(\sigma(\gamma))=\mu(\gamma)$ pour tout~$\gamma\in F$.

Dans ce paragraphe, nous remarquons qu'il existe un voisinage~$\mathcal{V}$ de \linebreak l'inclusion canonique dans $\Hom(\Gamma_0,G)$ tel que pour tout~$\sigma\in\mathcal{V}$ on ait \linebreak $\lambda(\sigma(\gamma))=\lambda(\gamma)$ et $\mu(\sigma(\gamma))=\mu(\gamma)$ pour tout~$\gamma\in\Gamma_0$, sans restriction à une partie finie~$F$.

\medskip

\begin{lem}\label{Reste un reseau cocompact sans torsion}
Soient $\kkk$ un corps local ultramétrique, $G$ l'ensemble des \linebreak $\kkk$-points d'un $\kkk$-groupe algébrique semi-simple connexe de $\kkk$-rang~un et $\Gamma_0$ un réseau cocompact sans torsion de~$G$.
Notons~$X$ l'arbre de Bruhat-Tits de~$G$ et $x_0$ le point de~$X$ donné par~\emph{(\ref{mu distance})}.
Pour tout domaine fondamental connexe~$\D$ de~$X$ pour l'action de~$\Gamma_0$, contenant~$x_0$, il existe un voisinage $\mathcal{V}$ de l'inclusion canonique dans $\Hom(\Gamma_0,G)$ tel que pour tout~$\sigma\in\mathcal{V}$,
\begin{itemize}
	\item le morphisme $\sigma$ soit injectif ;
	\item le groupe $\sigma(\Gamma_0)$ soit un réseau cocompact sans torsion de~$G$ ;
	\item l'ensemble~$\D$ soit un domaine fondamental de~$X$ pour l'action de~$\sigma(\Gamma_0)$ ;
	\item on ait $\lambda(\sigma(\gamma))=\lambda(\gamma)$ et $\mu(\sigma(\gamma))=\mu(\gamma)$ pour tout $\gamma\in\Gamma_0$.
\end{itemize}
\end{lem}

\begin{dem}
Soit~$\D$ un domaine fondamental connexe de~$X$ pour l'action de~$\Gamma_0$, contenant~$x_0$.
L'ensemble~$\mathcal{F}$ des éléments~$\gamma\in\Gamma_0$ tels que \linebreak $\gamma\cdot\D\cap\D\neq\emptyset$ est fini, donc
$$\D_1 = \bigcup_{\gamma\in\mathcal{F}} \gamma\cdot\D$$
est un compact de~$X$.
Le fixateur (point par point) de~$\D_1$ dans~$G$ est un voisinage de~$1$ dans~$G$, en tant qu'intersection des fixateurs des extrémités des arêtes de~$X$ rencontrant~$\D_1$.
Par conséquent, l'ensemble
$$\mathcal{V} = \big\{ \sigma\in\Hom(\Gamma_0,G),\quad \sigma(\gamma)\in\gamma\cdot K_{\D_1}\quad \forall\gamma\in\mathcal{F}\big\} $$
est un voisinage de l'inclusion canonique dans~$\Hom(\Gamma_0,G)$.

Pour voir que $\mathcal{V}$ vérifie les propriétés du lemme~\ref{Reste un reseau cocompact sans torsion}, il suffit de montrer que pour tout $\sigma\in\mathcal{V}$ il existe une isométrie bijective $\sigma$-équivariante $f_{\sigma} : X\rightarrow X$ fixant~$x_0$.
En effet, si une telle isométrie~$f_{\sigma}$ existe, alors $\sigma$ est injectif et $\sigma(\Gamma_0)$ est un réseau cocompact sans torsion de~$G$, admettant~$f_{\sigma}(\D)$ comme domaine fondamental dans~$X$.
Pour tout~$\gamma\in\Gamma_0$, l'image par~$f_{\sigma}$ de l'axe de translation~$\mathcal{A}_{\gamma}$ de~$\gamma$ est l'axe de translation~$\mathcal{A}_{\sigma(\gamma)}$ de~$\sigma(\gamma)$, et $\lambda(\sigma(\gamma))=\lambda(\gamma)$.
Enfin, comme~$f_{\sigma}$ est une isométrie $\sigma$-équivariante fixant~$x_0$, on a
$$\mu(\sigma(\gamma)) = d(x_0,\sigma(\gamma)\cdot x_0) = d\big(f_{\sigma}(x_0),f_{\sigma}(\gamma\cdot x_0)\big) = d(x_0,\gamma\cdot x_0) = \mu(\gamma)$$
pour tout $\gamma\in\Gamma_0$ d'après~(\ref{mu distance}).

Fixons donc un morphisme~$\sigma\in\mathcal{V}$ et établissons l'existence d'une isométrie bijective $\sigma$-équivariante $f_{\sigma} : X\rightarrow X$ fixant~$x_0$.
Pour tout~$x\in\D$ et tout~$\gamma\in\Gamma_0$ tels que~$\gamma\cdot x\in\D$, on a~$\sigma(\gamma)\cdot\nolinebreak x=\gamma\cdot x$.
Par conséquent, on peut définir une application $f_{\sigma} : X\rightarrow X$ en posant $f_{\sigma}(\gamma\cdot\nolinebreak x)=\sigma(\gamma)\cdot x$ pour tout~$\gamma\in\Gamma_0$ et tout~$x\in\D$.
Par construction, $f_{\sigma}$ est $\sigma$-équivariante et fixe~$x_0$.
Montrons que c'est une isométrie bijective.
Comme la restriction de~$f_{\sigma}$ à~$\D_1$ est l'identité, comme $f_{\sigma}$ est $\sigma$-équivariante et comme $\Gamma_0$ et~$\sigma(\Gamma_0)$ agissent sur~$X$ par isométries, la restriction de~$f_{\sigma}$ à tout translaté~$\gamma\cdot\D_1$, où~$\gamma\in\Gamma_0$, est une isométrie.
Comme les translatés par~$\Gamma_0$ de l'intérieur de~$\D_1$ recouvrent~$X$, on en déduit que $f_{\sigma}$ est une isométrie locale.
Comme $X$ est un arbre, c'est une isométrie globale.
Enfin, $f_{\sigma}$ est surjective.
En effet, la réunion~$X_{\D}$ des translatés par~$\sigma(\Gamma_0)$ de~$\D$ est ouverte dans~$X$, en tant que réunion des translatés par~$\sigma(\Gamma_0)$ de l'intérieur de~$\D_1$.
Elle est également fermée car toute suite de points de~$X_{\D}$ qui converge dans~$X$ est contenue à partir d'un certain rang dans une boule de~$X$ de diamètre~$\leq r$, où~$r>0$ désigne la distance de~$\D$ à~$X\smallsetminus\D_1$, donc dans un translaté par~$\sigma(\Gamma_0)$ de~$\D_1$ (car $\sigma(\Gamma_0)$ agit sur~$X$ par isométries).
Par connexité de~$X$, on a~$X_{\D}=X$, et donc
$$f_{\sigma}(X) = \bigcup_{\gamma\in\Gamma_0} f_{\sigma}(\gamma\cdot\D) = \bigcup_{\gamma\in\Gamma_0} \sigma(\gamma)\cdot f_{\sigma}(\D) = X_{\D} = X.$$
Ainsi, $f_{\sigma}:X\rightarrow X$ est une isométrie bijective $\sigma$-équivariante fixant~$x_0$, ce qui termine la démonstration.
\end{dem}

\medskip

Notons que l'existence de déformations non triviales est spécifique au rang~un.
En effet, si $\mathrm{rang}_{\kkk}(\mathbf{G})\geq 2$ et si par exemple $G$ est adjoint sans facteur compact et $\Gamma_0$ irréductible, alors tout morphisme $\varphi : \Gamma_0\rightarrow G$ tel que $\varphi(\Gamma_0)$ soit un réseau cocompact de~$G$ se prolonge en un automorphisme continu de~$G$ par le théorème de super-rigidité de Margulis (\cite{mar},~th.~5.6).

%%%%%%%%%%%%%%%%%%%%%%%%%
\subsection{Démonstration du théorème \ref{ouvert}}

Le théorème \ref{ouvert} est une conséquence du théorème~\ref{theoreme quotients compacts}, du théorème~\ref{condition d'admissibilite} et du lemme~\ref{Reste un reseau cocompact sans torsion}.
En effet, soit~$\Gamma$ un sous-groupe discret de type fini sans torsion de~$G\times G$ agissant proprement et cocompactement sur~$(G\times G)/\Delta_G$.
D'après le théorème~\ref{theoreme quotients compacts}, il existe un réseau cocompact sans torsion~$\Gamma_0$ de~$G$ et un morphisme admissible~$\rho : \Gamma_0\rightarrow G$ tels que
$$\Gamma = \big\{ (\gamma,\rho(\gamma)),\ \gamma\in\Gamma_0\big\} ,$$
à la permutation près des deux facteurs de~$G\times G$.
Soit~$X$ l'arbre de Bruhat-Tits de~$G$, soit~$N$ le nombre de sommets du graphe fini~$\Gamma_0\backslash X$ et soit~$F$ l'ensemble des éléments~$\gamma\in\Gamma_0\smallsetminus\{ 1\} $ tels que $\mu(\gamma)\leq 4N$.
Choisissons un domaine fondamental connexe~$\D$ de~$X$ pour l'action de~$\Gamma_0$ : on peut prendre par exemple le domaine de Dirichlet
$$\D = \big\{ x\in X,\quad d(x,x_0)\leq d(x,\gamma\cdot x_0)\quad \forall\gamma\in\Gamma_0\big\} .$$
Soit~$\mathcal{V}$ le voisinage de l'inclusion canonique dans~$\Hom(\Gamma_0,G)$ donné par le lemme~\ref{Reste un reseau cocompact sans torsion}.
D'après le théorème~\ref{condition d'admissibilite}, l'ensemble
$$\mathcal{W} = \big\{ \tau\in\Hom(\Gamma_0,G),\quad \lambda(\tau(\gamma))<\lambda(\gamma)\quad \forall\gamma\in F\big\} $$
contient~$\rho$ ; c'est donc un voisinage de~$\rho$ dans~$\Hom(\Gamma_0,G)$, et $\mathcal{V}\times\mathcal{W}$ est un voisinage de l'inclusion canonique dans $\Hom(\Gamma,G\times G)$.
Soit $\varphi=(\sigma,\tau)\in\nolinebreak\mathcal{V}\times\nolinebreak\mathcal{W}$.
Le groupe
$$\varphi(\Gamma) = \big\{ (\sigma(\gamma),\tau(\gamma)),\ \gamma\in\Gamma_0\big\} $$
est un graphe par injectivité de $\sigma$, et $\sigma(\Gamma_0)$ est un réseau cocompact sans torsion de~$G$.
Par définition de~$\mathcal{V}$ et de~$\mathcal{W}$, le graphe fini~$\sigma(\Gamma_0)\backslash X$ possède $N$ sommets, l'ensemble~$F$ est formé des éléments~$\gamma\in\Gamma_0\smallsetminus\{ 1\} $ tels que $\mu(\sigma(\gamma))\leq 4N$, et l'on a $\lambda(\tau(\gamma))<\lambda(\sigma(\gamma))$ pour tout~$\gamma\in F$.
D'après le théorème~\ref{condition d'admissibilite}, le morphisme $\tau\circ\sigma^{-1} : \sigma(\Gamma_0)\rightarrow G$ est admissible.
D'après le théorème~\ref{theoreme quotients compacts}, le groupe~$\varphi(\Gamma)$ agit librement, proprement et cocompactement sur~$(G\times G)/\Delta_G$.
\hfill\qedsymbol

%%%%%%%%%%%%%%%%%%%%%%%%%%%%%%%%%%%%%%%%%%%%%%%%%%%%%%%%%%%%
\section{Lien avec l'outre-espace}\label{outre-espace}

Dans le cas particulier où~$\rho$ est injectif d'image discrète et cocompacte, la proposition~\ref{max sur un ensemble fini} implique l'existence et l'équivalence entre deux définitions différentes d'une même ``distance asymétrique'' sur l'outre-espace, comme nous le précisons dans cette partie.

Cette distance asymétrique sur l'outre-espace est un analogue de la distance asymétrique de Thurston sur l'espace de Teichmüller, introduite et étudiée en détail dans~\cite{thu} : les classes d'équivalence de structures hyperboliques complètes sur une surface compacte~$S$ de type fini et de caractéristique d'Euler~$\chi(S)<0$,
munies d'un marquage par le groupe fondamental de~$S$, sont remplacées par les classes d'équivalence de graphes métriques connexes munis d'un marquage par un groupe libre de type fini fixé (\textit{cf.} paragraphe~\ref{lien avec Teichmuller}).

Sur l'outre-espace, cette distance asymétrique a d'abord été étudiée par T.~White dans un texte non publié, puis récemment par S.~Francaviglia et A.~Martino~\cite{fm} qui se sont intéressés à une symétrisation de cette distance et à la géométrie correspondante sur l'outre-espace.

%%%%%%%%%%%%%%%%%%%%%%%%%
\subsection{Rappels}

Soient~$n\geq 2$ un entier et $\mathcal{R}_n$ un bouquet de $n$ cercles d'intersection~$\{ c\} $ ; le groupe fondamental~$\pi_1(\mathcal{R}_n,c)$ est un groupe libre à~$n$ générateurs que nous noterons~$\F_n$.
Appelons \textit{graphe normalisé marqué de rang}~$n$ tout couple~$(Y,\varphi)$ où $Y$ est un graphe métrique connexe fini de valence~$\geq 2$ dont la somme des longueurs des arêtes vaut un, et où $\varphi : \mathcal{R}_n\rightarrow Y$ est une équivalence d'homotopie.
L'équivalence d'homotopie~$\varphi$ induit un isomorphisme de groupes
$$\varphi_{\ast} : \F_n \longrightarrow \pi_1\big(Y,\varphi(c)\big).$$
Pour tous graphes normalisés~$(Y,\varphi)$ et~$(Y',\varphi')$ marqués de rang~$n$, nous dirons qu'une application continue $\overline{f} : Y\rightarrow Y'$ \textit{respecte les marquages~$\varphi$ et~$\varphi'$} si $\overline{f}\circ\varphi$ est homotope à~$\varphi'$.
Notons~$(Y,\varphi)\sim(Y',\varphi')$ s'il existe
une isométrie bijective $i : Y\rightarrow Y'$ respectant les marquages~$\varphi$ et~$\varphi'$.
Ceci définit une relation d'équivalence~$\sim$ sur l'ensemble des graphes normalisés marqués de rang~$n$.
On appelle \textit{outre-espace} de rang~$n$ l'ensemble des classes d'équivalence pour cette relation.
Cet ensemble a été introduit par M.~Culler et K.~Vogtmann dans l'article fondateur~\cite{cv} ; nous le notons ici~$\OS$.
Pour tout graphe normalisé~$(Y,\varphi)$ marqué de rang~$n$, nous notons~$[Y,\varphi]$ la classe de~$(Y,\varphi)$ dans~$\OS$.

%%%%%%%%%%%%%%%%%%%%%%%%%
\subsection{Une distance asymétrique sur l'outre-espace}\label{distance asymetrique}

Soient~$(Y,\varphi)$ et $(Y',\varphi')$ deux graphes normalisés marqués de rang~$n$.
Notons~$X$ (resp.~$X'$) un revêtement universel de~$Y$ (resp. de~$Y'$) : c'est un arbre réel simplicial de valence~$\geq 2$.
Le groupe fondamental~$\pi_1(Y,\varphi(c))$ (resp.~$\pi_1(Y',\varphi'(c))$) agit librement, proprement et cocompactement sur~$X$ (resp. sur~$X'$), par isométries.
Posons~$\Gamma_0=\pi_1(Y,\varphi(c))$ et
$$\rho = \varphi'_{\ast}\circ\varphi_{\ast}^{-1}\ :\ \Gamma_0\longrightarrow\Isom(X').$$
Toute application continue $\overline{f} : Y\rightarrow Y'$ respectant les marquages~$\varphi$ et~$\varphi'$ se relève en une application continue $\rho$-équivariante $f : X\rightarrow X'$, au sens de la proposition~\ref{max sur un ensemble fini}.
Réciproquement, toute application continue $\rho$-équivariante $f : X\rightarrow X'$ induit une application continue $\overline{f} : Y\rightarrow Y'$ respectant les marquages~$\varphi$ et~$\varphi'$.
Par construction, l'application~$f$ est lipschitzienne si et seulement si~$\overline{f}$ l'est, et dans ce cas les deux constantes de Lipschitz sont les mêmes.

D'après la proposition~\ref{max sur un ensemble fini}, il existe une application $\overline{f} : Y\rightarrow Y'$ affine par morceaux respectant les marquages~$\varphi$ et~$\varphi'$, et la borne inférieure des constantes de Lipschitz de telles applications est atteinte.
Cette borne inférieure ne dépend que des classes $[Y,\varphi]$ et~$[Y',\varphi']$ dans~$\OS$ ; notons $L([Y,\varphi],[Y',\varphi'])$ son logarithme.

\begin{lem}\label{distance sur l'outre-espace, applications lipschitziennes}
L'application $L : \OS\times\OS\rightarrow\R$ est une distance asymétrique sur~$\OS$, au sens où
\begin{enumerate}
	\item pour tous~$[Y,\varphi],[Y',\varphi']\in\OS$ on a
	$$L([Y,\varphi],[Y',\varphi'])\geq 0,$$
	avec égalité si et seulement si~$[Y,\varphi]=[Y',\varphi']$ ;
	\item pour tous $[Y,\varphi],[Y',\varphi'],[Y'',\varphi'']\in\OS$ on a
	$$L([Y,\varphi],[Y'',\varphi''])\,\leq\,L([Y,\varphi],[Y',\varphi'])\,+\,L([Y',\varphi'],[Y'',\varphi'']).$$
\end{enumerate}
\end{lem}

\medskip

En général on a $L([Y,\varphi],[Y',\varphi'])\neq L([Y',\varphi'],[Y,\varphi])$, comme nous le verrons ci-dessous.

\medskip

\noindent
\textbf{Démonstration du lemme~\ref{distance sur l'outre-espace, applications lipschitziennes}.}
Le point~(2) est clair.
Prouvons le point~(1).
Soient $(Y,\varphi)$ et~$(Y',\varphi')$ deux graphes normalisés marqués de rang~$n$ et $\overline{f} : Y\rightarrow Y'$ une application affine par morceaux respectant les marquages~$\varphi$ et~$\varphi'$, de constante de Lipschitz $C=\exp L([X,\varphi],[X',\varphi'])$ minimale.

Remarquons que $\overline{f}$ est surjective.
En effet, notons comme ci-dessus $X$ (resp.~$X'$) un revêtement universel de~$Y$ (resp. de~$Y'$) et soit $f : X\rightarrow X'$ un relevé de~$\overline{f}$ ; montrons que~$f$ est surjective.
Fixons un sommet~$s$ de~$X$.
Pour tout~$\gamma\in\nolinebreak\Gamma_0$, l'image par~$f$ du segment géodésique~$[s,\gamma\cdot s]$ de~$X$ contient le segment géodésique~$[f(s),\rho(\gamma)\cdot f(s)]$ de~$X'$.
Il suffit donc de montrer que les segments géodésiques~$[f(s),\rho(\gamma)\cdot f(s)]$, où $\gamma\in\Gamma_0$, recouvrent~$X'$.
Mais s'il existait un point~$x'\in X'$ n'appartenant à aucun segment géodésique~$[f(s),\rho(\gamma)\cdot f(s)]$, alors toute composante connexe de~$X'\smallsetminus\{ x'\} $ ne contenant pas~$f(s)$ serait un sous-arbre infini de~$X$ ne rencontrant pas l'orbite~$\rho(\Gamma_0)\cdot f(s)$, ce qui contredirait le fait que~$\rho(\Gamma_0)$ agit cocompactement sur~$X'$.

Notons~$E$ (resp.~$E'$) l'ensemble des arêtes de~$Y$ (resp. de~$Y'$).
Pour toute arête~$e\in E$ (resp.~$e'\in E'$), notons~$\ell(e)\in [0,1]$ (resp.~$\ell'(e')\in [0,1]$) sa longueur dans~$Y$ (resp. dans~$Y'$).
Notons enfin~$C_e$ la constante de Lipschitz de la restriction de~$\overline{f}$ à~$e$ pour tout~$e\in E$.
Comme $\overline{f}$ est surjective on a
\begin{equation}\label{inegalites constantes de Lipschitz}
1\, =\, \sum_{e'\in E} \ell'(e')\, \leq\, \sum_{e\in E} \ell\big(\overline{f}(e)\big)\, \leq\, \sum_{e\in E} C_e\,\ell(e)\, \leq\, C \sum_{e\in E} \ell(e)\, =\, C.
\end{equation}
On en déduit
\begin{equation}\label{L positif}
L([X,\varphi],[X',\varphi']) = \log C \geq 0.
\end{equation}
De plus, si~$L([X,\varphi],[X',\varphi'])=0$, alors toutes les inégalités dans~(\ref{inegalites constantes de Lipschitz}) sont des égalités ; on en déduit aisément que~$\overline{f}$ est une isométrie bijective, et donc que~$[X,\varphi]=[X',\varphi']$.
Enfin, en considérant l'identité de~$X$, on voit que~$L([X,\varphi],[X,\varphi])\leq 0$, donc~$L([X,\varphi],[X,\varphi])=0$ par~(\ref{L positif}).
\hfill\qedsymbol

\bigskip

Montrons l'asymétrie sur un exemple pour~$n=2$.
Comme précédemment, soit~$\mathcal{R}_2$ un bouquet de deux cercles~$\mathcal{C}$ et~$\mathcal{C}'$ d'intersection~$\{ c\} $.
Pour~$i\in\{ 1,2\} $, soit~$Y_i$ le graphe donné par la figure~5, où les deux boucles sont de longueur $a_i\in]0,1[$ et où l'arête du milieu est de longueur $1-2a_i$.
Soit $\varphi_i : \mathcal{R}_2\rightarrow Y_i$ une équivalence d'homotopie envoyant~$c$ sur~$y_i$, établissant un homéomorphisme entre~$\mathcal{C}$ et le cercle de gauche de~$Y_i$ et envoyant~$\mathcal{C}'$ sur l'union du segment transverse et du cercle de droite de~$Y_i$.
Pour~$a_1\leq a_2$, un calcul donne
$$L\big([Y_1,\varphi_1],[Y_2,\varphi_2]\big)=\log\Big(\frac{a_2}{a_1}\Big) \quad\mathrm{et}\quad L\big([Y_2,\varphi_2],[Y_1,\varphi_1]\big)=\log\bigg(\frac{1-a_1}{1-a_2}\bigg).$$
Ces réels sont différents dès que $a_2\notin\{ a_1,1-a_1\} $.

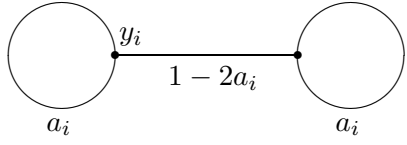
\begin{figure}\label{Graphe}
\begin{center}\setlength{\unitlength}{1mm}
\begin{picture}(50,20)(-25,-10)
\put(-19,0){\circle{14}}
\put(-12,0){\circle*{1}}
\put(-12,0){\line(1,0){24}}
\put(12,0){\circle*{1}}
\put(19,0){\circle{14}}

\put(-21,-10){$a_i$}
\put(-5,-4){$1-2a_i$}
\put(17,-10){$a_i$}
\put(-11.5,2){$y_i$}
\end{picture}
\end{center}
\caption{\textit{Le graphe $Y_i$}}
\end{figure}

%%%%%%%%%%%%%%%%%%%%%%%%%
\subsection{Quotients de longueurs de lacets}\label{quotients lacets}

Pour tout graphe normalisé~$(Y,\varphi)$ marqué de rang~$n$ et tout~$\gamma\in\F_n\smallsetminus\{ 1\} $, la longueur minimale d'un lacet dans la classe d'homotopie libre de~$\varphi_{\ast}(\gamma)\in\pi_1(Y,\varphi(c))$ ne dépend que de la classe~$[Y,\varphi]$ de~$(Y,\varphi)$ dans~$\OS$ ; notons-la~$\mathrm{long}_{[Y,\varphi]}(\gamma)$.
Si~$X$ désigne comme précédemment un revêtement universel de~$Y$, alors~$\mathrm{long}_{[Y,\varphi]}(\gamma) = \lambda(\varphi_{\ast}(\gamma))$ est la longueur de translation de~$\varphi_{\ast}(\gamma)$ vu comme isométrie de~$X$.

Pour tous $[Y,\varphi],[Y',\varphi']\in\OS$, tout~$\gamma\in\F_n\smallsetminus\{ 1\} $ et toute application continue $\overline{f} : Y\rightarrow Y'$ respectant les marquages~$\varphi$ et~$\varphi'$, lipschitzienne de constante~$C$ minimale, on a
\begin{equation}\label{inegalite lambda C}
\frac{\mathrm{long}_{[Y',\varphi']}(\gamma)}{\mathrm{long}_{[Y,\varphi]}(\gamma)}\ \leq\ C.
\end{equation}
Posons
$$K\big([Y,\varphi],[Y',\varphi']\big) = \log \sup_{\gamma\in\F_n\smallsetminus\{ 1\} } \bigg(\frac{\mathrm{long}_{[Y',\varphi']}(\gamma)}{\mathrm{long}_{[Y,\varphi]}(\gamma)}\bigg).$$
Les inégalités~(\ref{inegalite lambda C}) impliquent
$$K\big([Y,\varphi],[Y',\varphi']\big)\ \leq\ L\big([Y,\varphi],[Y',\varphi']\big).$$
D'après la proposition~\ref{max sur un ensemble fini}, cette inégalité est en fait une égalité et la borne supérieure $K([Y,\varphi],[Y',\varphi'])$ est atteinte.

\medskip

\begin{coro}\label{K=L}
\emph{(T. White,} cf. \emph{\cite{fm}, prop.~3.11)}\\
Pour tous~$[Y,\varphi],[Y',\varphi']\in\nolinebreak\OS$ il existe un élément~$\gamma\in\F_n\smallsetminus\{ 1\} $ tel que
$$e^{K([Y,\varphi],[Y',\varphi'])}\ =\ \frac{\mathrm{long}_{[Y',\varphi']}(\gamma)}{\mathrm{long}_{[Y,\varphi]}(\gamma)}\ =\ \frac{\lambda\big(\varphi'_{\ast}(\gamma)\big)}{\lambda\big(\varphi_{\ast}(\gamma)\big)}\ =\ e^{L([Y,\varphi],[Y',\varphi'])}\ \geq\ 1.$$
\end{coro}

En particulier on a~$K=L$.

%%%%%%%%%%%%%%%%%%%%%%%%%
\subsection{Lien avec la distance asymétrique de Thurston sur l'espace de\\ Teichmüller}\label{lien avec Teichmuller}

On a une analogie très forte entre la distance asymétrique $K=L$ sur l'outre-espace et la distance asymétrique de Thurston sur l'espace de Teichmüller.
Plus précisément, soit~$S$ une surface compacte de caractéristique d'Euler~$\chi(S)<0$, et soit~$\mathcal{T}(S)$ l'espace de Teichmüller de~$S$, défini comme l'ensemble des classes d'équivalence de structures hyperboliques complètes sur~$S$ pour la relation ``être tirée en arrière par un homéomorphisme de~$S$ homotope à l'identité''.

Pour toutes structures hyperboliques complètes~$g$ et~$h$ sur~$S$, la borne inférieure des constantes de Lipschitz d'homéomorphismes~$\overline{f} : (S,g)\rightarrow (S,h)$ homotopes à l'identité ne dépend que des classes de~$g$ et~$h$ dans~$\mathcal{T}(S)$.
Notons~$L([g],[h])$ le logarithme de cette borne inférieure.
D'après~\cite{thu}, prop.~2.1, l'application~$L$ est une distance asymétrique sur~$\mathcal{T}(S)$, au sens du lemme~\ref{distance sur l'outre-espace, applications lipschitziennes}.

Soit~$x_0\in S$ un point-base.
Pour toute structure hyperbolique complète~$g$ sur~$S$ et tout élément~$\gamma\in\pi_1(S,x_0)$ non trivial, notons~$\mathrm{long}_g(\gamma)>0$ la plus petite longueur d'un lacet dans la classe d'homotopie libre de~$\gamma$ ; elle est atteinte par l'unique géodésique fermée pour~$g$ dans la classe d'homotopie libre de~$\gamma$.
Ceci définit une application~$\mathrm{long}_g$ qui ne dépend que de la classe de~$g$ dans~$\mathcal{T}(S)$.
D'après~\cite{thu}, th.~3.1, l'application~$K : \mathcal{T}(S)\times\mathcal{T}(S)\rightarrow\R$ définie par
$$K([g],[h]) = \log \sup_{\gamma\in\pi_1(S,x_0)\smallsetminus\{ 1\} } \bigg(\frac{\mathrm{long}_h(\gamma)}{\mathrm{long}_g(\gamma)}\bigg)$$
pour tous~$[g],[h]\in\mathcal{T}(S)$ est une distance asymétrique sur~$\mathcal{T}(S)$.

Comme sur l'outre-espace, on a~$K=L$ d'après~\cite{thu}, th.~8.5.
En revanche, dans le cas de l'espace de Teichmüller la borne supérieure~$K([g],[h])$ n'est en général pas atteinte par une géodésique fermée simple, mais par une lamination géodésique mesurée.

%%%%%%%%%%%%%%%%%%%%%%%%%
\subsection{Absence de morphisme admissible injectif d'image discrète}\label{Image non discrete}

D'après le corollaire~\ref{K=L} et le lemme~\ref{distance sur l'outre-espace, applications lipschitziennes}, on a $K([Y,\varphi],[Y',\varphi'])\geq 0$ pour tous $[Y,\varphi],[Y',\varphi']\in\OS$, et en cas d'égalité on a $[Y,\varphi]=[Y',\varphi']$ donc $\lambda(\varphi'_{\ast}(\gamma))=\lambda(\varphi_{\ast}(\gamma))$ pour tout~$\gamma\in\F_n$.
Nous allons en déduire le résultat suivant, dont le corollaire~\ref{fidele et discrete} est un cas particulier.

\begin{coro}\label{fidele et discrete sur les arbres}
Soit $X$ un arbre réel simplicial, bipartite de valences $n_1\geq 2$ et $n_2\geq 2$, dont toutes les arêtes ont même longueur, et soit~$\Gamma_0$ un sous-groupe discret sans torsion de~$\Isom(X)$ tel que le graphe~$\Gamma_0\backslash X$ soit fini.
Il n'existe pas de morphisme de groupes admissible $\rho : \Gamma_0\rightarrow\nolinebreak\Isom(X)$ qui soit injectif d'image discrète.
\end{coro}

Soient~$x_0,x'_0\in X$.
Rappelons qu'un morphisme $\rho : \Gamma_0\rightarrow\Isom(X)$ est dit \textit{admissible} si pour tout~$R>0$ on a $d(x'_0,\rho(\gamma)\cdot x'_0)\leq d(x_0,\gamma\cdot x_0)-R$ pour presque tout~$\gamma\in\Gamma_0$ (\textit{cf.} paragraphe~\ref{Demonstration du theoreme sur l'admissibilite}).
Cette condition ne dépend pas du choix de~$x_0$ et~$x'_0$.

De même, la positivité de la distance asymétrique de Thurston sur l'espace de Teichmüller (avec la condition d'égalité) implique l'absence de morphisme admissible pour~$G=\PSL_2(\R)$ (\textit{cf.} \cite{sa2}, \S~4.1).
On n'a pas besoin que la borne supérieure~$K([g],[h])$ soit atteinte par une géodésique fermée simple.

Pour démontrer le corollaire~\ref{fidele et discrete sur les arbres} nous utilisons le lemme suivant.

\begin{lem}\label{nombre d'aretes}
Le nombre d'arêtes d'un graphe connexe fini dont le groupe fondamental est libre de rang $n\geq 2$ fixé est une fonction décroissante de la valence moyenne des sommets.
\end{lem}

\begin{dem}
Soit~$Y$ un graphe connexe fini dont le groupe fondamental est libre de rang~$n\geq 2$.
Si~$a$ désigne le nombre d'arêtes, $s$ le nombre de sommets et $v$ la valence moyenne des sommets de~$Y$, alors~$2a = vs$.
De plus, d'après~\cite{die}, th.~1.9.6, on a $n = a-s+1$, d'où
$$a = \frac{n-1}{1-\frac{2}{v}},$$
ce qui implique le lemme.
\end{dem}

\medskip

Nous pouvons à présent démontrer le corollaire~\ref{fidele et discrete sur les arbres}.

\medskip

\noindent
\textbf{Démonstration du corollaire~\ref{fidele et discrete sur les arbres}.}
Soit $\rho : \Gamma_0\rightarrow\Isom(X)$ un morphisme de groupes injectif d'image discrète.
Montrons que~$\rho$ n'est pas admissible.

Les groupes~$\Gamma_0$ et~$\rho(\Gamma_0)$ agissent tous deux librement et proprement sur~$X$.
Les graphes quotients~$\Gamma_0\backslash X$ et~$\rho(\Gamma_0)\backslash X$ sont connexes, bipartites de valences $n_1\geq 2$ et $n_2\geq 2$, et leurs groupes fondamentaux, qui s'identifient respectivement à~$\Gamma_0$ et~$\rho(\Gamma_0)$, sont libres de même rang.
Le graphe~$\Gamma_0\backslash X$ est fini mais~$\rho(\Gamma_0)\backslash X$ ne l'est pas forcément.

Comme~$\Gamma_0$ est sans torsion et comme~$\rho$ est injectif d'image discrète, l'élément~$\rho(\gamma)$ est hyperbolique pour tout~$\gamma\in\Gamma_0\smallsetminus\{ 1\} $.
Fixons un élément~$\gamma\in\Gamma_0\smallsetminus\{ 1\} $ et un sommet~$x'\in X$ de l'axe de translation~$\mathcal{A}_{\rho(\gamma)}$.
L'image~$Y'$ dans~$\rho(\Gamma_0)\backslash X$ de l'union des segments géodésiques~$[x',\rho(\gamma)\cdot x']$, où~$\gamma\in\Gamma_0$, est un sous-graphe fini connexe de~$\rho(\Gamma_0)\backslash X$ dont le groupe fondamental s'identifie encore à~$\rho(\Gamma_0)$.
Si~$\rho(\Gamma_0)\backslash X$ est fini, on voit facilement que~$Y'=\rho(\Gamma_0)\backslash X$.
Par construction, $Y'$ est de valence~$\geq 2$.
De plus, la valence moyenne des sommets de~$Y'$ est inférieure à~$(n_1+n_2)/2$ : en effet, si l'on note~$S_1$ (resp.~$S_2$) l'ensemble des sommets de~$\rho(\Gamma_0)\backslash X$ de valence~$n_1$ (resp.~$n_2$), alors $Y'\cap S_1$ et~$Y'\cap S_2$ ont même cardinal, et pour tout~$s\in Y'\cap S_i$, la valence de~$s$ dans~$Y'$ est inférieure à~$n_i$.

Notons~$m$ (resp.~$m'$) le nombre d'arêtes de~$\Gamma_0\backslash X$ (resp. de~$Y'$).
D'après le lemme~\ref{nombre d'aretes}, on a~$m\leq m'$.
Or, toutes les arêtes de~$X$ ont même longueur par hypothèse.
Si l'on note~$L>0$ cette longueur et si l'on munit~$\Gamma_0\backslash X$ (resp.~$Y'$) de la distance induite par celle de~$X$ divisée par~$mL$ (resp. par~$m'L$), alors la somme des longueurs des arêtes de~$\Gamma_0\backslash X$ (resp. de~$Y'$) vaut~un.
D'après le corollaire~\ref{K=L}, il existe un élément~$\gamma\in\Gamma_0$ non trivial tel que
$$\frac{1}{m'L}\lambda(\rho(\gamma))\ \geq\ \frac{1}{mL}\lambda(\gamma).$$
Comme~$m\leq m'$, on en déduit $\lambda(\rho(\gamma))\geq \lambda(\gamma)$.
D'après le théorème~\ref{condition d'admissibilite sur les arbres}, le morphisme~$\rho$ n'est pas admissible.
\hfill\qedsymbol

\bigskip

Le corollaire~\ref{fidele et discrete} s'obtient à partir du corollaire~\ref{fidele et discrete sur les arbres} en prenant pour~$X$ l'arbre de Bruhat-Tits de~$G$.

%%%%%%%%%%%%%%%%%%%%%%%%%%%%%%%%%%%%%%%%%%%%%%%%%%%%%%%%%%%%
\section{Existence de quotients compacts par $\Gamma$ Zariski-dense}\label{Zariski-dense}

Pour démontrer le corollaire~\ref{coro Zariski-dense}, nous utilisons le résultat suivant, dont l'analogue réel est dû à T.~Kobayashi (\cite{ko3}, th.~2.4).

\begin{lem}\label{lemme Kobayashi}
Soient $\kkk$ un corps local ultramétrique, $G$ l'ensemble des \linebreak $\kkk$-points d'un $\kkk$-groupe algébrique semi-simple connexe de $\kkk$-rang~un et $\Gamma_0$ un réseau cocompact sans torsion de~$G$.
Soient~$X$ l'arbre de Bruhat-Tits de~$G$ et~$\D$ un domaine fondamental connexe de~$X$ pour l'action de~$\Gamma_0$, contenant dans son intérieur le point~$x_0$ donné par~\emph{(\ref{mu distance})}.
Posons
$$\mathcal{F} = \big\{ \gamma\in\Gamma_0\smallsetminus\{ 1\} ,\quad \gamma\cdot\D\cap\D\neq\emptyset\big\} $$
et notons~$\delta>0$ la distance de~$\D$ au complémentaire de~$\bigcup_{\gamma\in\mathcal{F}} \gamma\cdot\D$ dans~$X$.
Alors tout morphisme de groupes $\rho : \Gamma_0\rightarrow G$ vérifiant~$\mu(\rho(\gamma))<\delta$ pour tout~$\gamma\in\mathcal{F}$ est admissible.
\end{lem}

Le raisonnement de Kobayashi se transpose au cas ultramétrique en remplaçant simplement l'espace symétrique~$G/K$ par l'arbre de Bruhat-Tits de~$G$ ; nous le reproduisons ici pour la commodité du lecteur.

\medskip

\noindent
\textbf{Démonstration du lemme~\ref{lemme Kobayashi}.}
Soit~$\gamma\in\Gamma_0$.
Notons~$n$ la partie entière de~$d(x_0,\gamma\cdot x_0)/\delta$ et choisissons une suite~$(x_i)_{i=1,\ldots,n+1}$ de points du segment géodésique~$[x_0,\gamma\cdot x_0]$ telle que~$x_{n+1}=\gamma\cdot x_0$ et telle que pour tout~$0\leq i\leq n$ on ait~$d(x_i,x_{i+1})<\delta$.
Par récurrence, on construit une suite~$(\gamma_i)_{i=0,\ldots,n}$ d'éléments de~$\Gamma_0$ telle que $x_{i+1}\in\nolinebreak\gamma_0\ldots\gamma_i\cdot\D$ pour tout~$0\leq i\leq n$ ; par définition de~$\delta$ et comme $d(x_i,x_{i+1})<\delta$, on a~$\gamma_i\in\mathcal{F}\cup\{ 1\} $ pour tout~$i$.
D'autre part, on a $x_{n+1}=\gamma\cdot x_0\in\gamma_0\ldots\gamma_n\cdot\D$, donc $x_0$ appartient à la fois à l'intérieur de~$\D$ et à $(\gamma^{-1}\gamma_0\ldots\gamma_n)\cdot\D$.
Comme $\D$ est un domaine fondamental de $X$ pour~$\Gamma_0$, on~a $\gamma^{-1}\gamma_0\ldots\gamma_n=1$, \textit{i.e.} $\gamma=\gamma_0\ldots\gamma_n$.
Si $\ell_{\mathcal{F}} : \Gamma_0\rightarrow\N$ désigne la longueur des mots associée à~$\mathcal{F}$, on a
\begin{equation}\label{longueur des mots}
\ell_{\mathcal{F}}(\gamma)\ \leq\ n+1\ \leq\ \frac{d(x_0,\gamma\cdot x_0)}{\delta}+1\ =\ \frac{\mu(\gamma)}{\delta}+1.
\end{equation}
Soit $\rho\in\Hom(\Gamma_0,G)$ tel que
$$M := \max\big\{ \mu(\rho(\gamma)),\ \gamma\in\mathcal{F}\big\}\ <\ \delta.$$
Pour tout $\gamma\in\Gamma_0$ on a, d'après~(\ref{inegalite triangulaire pour mu}) et~(\ref{longueur des mots}),
$$\mu(\rho(\gamma))\ \leq\ M\cdot\ell_{\mathcal{F}}(\gamma)\ \leq\ \frac{M}{\delta} \mu(\gamma) + M.$$
Pour tout~$R>0$ on a $\frac{M}{\delta} \mu(\gamma) + M > \mu(\gamma) - R$ si et seulement si $\mu(\gamma) < \frac{M+R}{1-M/\delta}$. L'ensemble des éléments $\gamma\in\Gamma_0$ vérifiant cette inégalité est fini car $\Gamma_0$ est discret dans $G$ et l'application $\mu$ est propre.
Ainsi, $\rho$ est admissible.
\hfill\qedsymbol

\bigskip

Nous pouvons à présent démontrer la proposition~\ref{coro Zariski-dense}.

\medskip

\noindent
\textbf{Démonstration de la proposition~\ref{coro Zariski-dense}.}
Notons~$L>0$ la longueur commune des arêtes de~$X$.
En s'inspirant de la démonstration du théorème~2.1 de~\cite{lub}, par exemple, on construit facilement un réseau cocompact sans torsion~$\Gamma_0$ de~$G$ admettant un domaine fondamental connexe~$\D$ dans~$X$ qui contient dans son intérieur le point~$x_0$ donné par~(\ref{mu distance}) et tel que, en posant
$$\mathcal{F} = \big\{ \gamma\in\Gamma_0\smallsetminus\{ 1\} ,\quad \gamma\cdot\D\cap\D\neq\emptyset\big\} ,$$
la distance~$\delta$ de~$\D$ au complémentaire de~$\bigcup_{\gamma\in\mathcal{F}} \gamma\cdot\D$ dans~$X$ soit supérieure à~$2L$.
Soit~$\mathcal{F}'$ une partie de~$\mathcal{F}$ telle que~$\mathcal{F}$ soit l'union disjointe de~$\mathcal{F}'$ et de~${\mathcal{F}'}^{-1}$.
Le groupe~$\Gamma_0$ est libre, librement engendré par~$\mathcal{F}'$, donc tout morphisme de~$\Gamma_0$ dans~$G$ est entièrement déterminé par son image sur~$\mathcal{F}'$.

Soient~$\gamma_1\neq\gamma_2$ deux éléments de~$\mathcal{F}'$.
Rappelons qu'un élément de~$G$ est dit \textit{régulier} si la composante neutre de son centralisateur dans~$\mathbf{G}$ est un tore maximal de~$\mathbf{G}$.
L'ensemble des éléments réguliers de~$G$ contient un ouvert de Zariski de~$G$ (\cite{bor}, th.~12.3).
Par conséquent, si l'on note~$\mathcal{V}$ le voisinage de l'inclusion canonique dans~$\Hom(\Gamma_0,G)$ donné par le lemme~\ref{Reste un reseau cocompact sans torsion}, alors $\mathcal{V}\cdot\gamma_1$ contient un élément régulier~$\gamma'_1$.
Comme~$\delta\geq 2L$, l'ensemble des éléments~$g\in\nolinebreak G$ tels que~$0<\lambda(g)\leq\mu(g)<\delta$ est un ouvert non vide de~$G$ ; il contient donc un élément régulier~$\gamma''_1$.
Par un résultat de J.~Tits (\cite{ti2}, prop.~4.4), la réunion de tous les sous-groupes stricts de $G\times G$ qui contiennent~$(\gamma'_1,\gamma''_1)$ et qui sont Zariski-fermés et Zariski-connexes est incluse dans un fermé de Zariski strict~$F_1$ de $G\times G$.
Comme précédemment, l'ouvert de Zariski $(G\times G)\smallsetminus F_1$ contient un élément régulier~$(\gamma'_2,\gamma''_2)$ tel que~$\gamma'_2\in\mathcal{V}\cdot\gamma_2$ et~$0<\lambda(\gamma''_2)\leq\mu(\gamma''_2)<\delta$.
Par construction, le groupe engendré par $(\gamma'_1,\gamma''_1)$ et $(\gamma'_2,\gamma''_2)$ est Zariski-dense dans~$G\times G$, et sa projection sur chacun des facteurs de~$G\times G$ est non bornée.

Soit $\sigma : \Gamma_0\rightarrow G$ le morphisme de groupes défini par $\sigma(\gamma_i)=\gamma'_i$ pour $i\in\{ 1,2\} $ et $\sigma(\gamma)=\gamma$ pour~$\gamma\in\mathcal{F}'\smallsetminus\{ \gamma_1,\gamma_2\} $.
Par construction, on a $\sigma\in\mathcal{V}$.
Pour tout~$\gamma\in\mathcal{F}'\smallsetminus\{ \gamma_1,\gamma_2\} $, choisissons un élément $g_{\gamma}\in G$ tel que~$\mu(g_{\gamma})<\delta$.
Soit $\rho : \Gamma_0\rightarrow G$ le morphisme de groupes défini par $\rho(\gamma_i)=\gamma''_i$ pour~$i\in\{ 1,2\} $ et $\rho(\gamma)=g_{\gamma}$ pour~$\gamma\in\mathcal{F}'\smallsetminus\{ \gamma_1,\gamma_2\} $.
Le groupe
$$\Gamma = \big\{ (\sigma(\gamma),\rho(\gamma)),\ \gamma\in\Gamma_0\big\} $$
est Zariski-dense dans~$G\times G$ et sa projection sur chacun des facteurs de~$G\times G$ est non bornée.
Comme $\sigma\in\mathcal{V}$, le morphisme~$\sigma$ est injectif, le groupe~$\sigma(\Gamma_0)$ est un réseau cocompact sans torsion de~$G$ de domaine fondamental~$\D$ dans~$X$ et $\delta$ est la distance de~$\D$ au complémentaire de~$\bigcup_{\gamma\in\mathcal{F}} \sigma(\gamma)\cdot\D$ dans~$X$.
Par construction on a $\mu(\rho(\gamma))<\delta$ pour tout~$\gamma\in\mathcal{F}'$, donc le morphisme \linebreak $\rho\circ\sigma^{-1} : \sigma(\Gamma_0)\rightarrow G$ est admissible d'après le lemme~\ref{lemme Kobayashi} et l'inégalité~(\ref{mu de l'inverse}).
Le théorème~\ref{theoreme quotients compacts} assure que le groupe~$\Gamma$ agit librement, proprement et cocompactement sur $(G\times G)/\Delta_G$, ce qui termine la démonstration de la proposition~\ref{coro Zariski-dense}.
\hfill \qedsymbol

\vspace{0.5cm}
%%%%%%%%%%%%%%%%%%%%%%%%%%%%%%%%%%%%%%%%%%%%%%%%%%%%%%%%%%%%

\vspace{0.5cm}

{\small \textsc{D\'epartement de Math\'ematiques,
B\^atiment~425,
Facult\'e des Sciences d'Orsay,
Universit\'e Paris-Sud 11,
91405 Orsay Cedex,
France}

\textit{Adresse \'electronique :}
\texttt{fanny.kassel@math.u-psud.fr}}

\end{document}